\documentclass[twoside, letterpaper, 12pt]{article}

\DeclareMathSizes{10}{10}{10}{10}
\usepackage{anysize}
\usepackage{fancyhdr}
\marginsize{2cm}{2cm}{0.5cm}{0.75cm}
\usepackage{setspace}
\usepackage[utf8]{inputenc} % allow utf-8 input
\usepackage[T1]{fontenc}    % use 8-bit T1 fonts
\usepackage[colorlinks]{hyperref}       % hyperlinks
\usepackage{url}            % simple URL typesetting
\usepackage{booktabs}       % professional-quality tables
\usepackage{amsfonts}       % blackboard math symbols
\usepackage{nicefrac}       % compact symbols for 1/2, etc.
\usepackage{microtype}      % microtypography

\usepackage{caption}
\usepackage{subcaption}
\usepackage{graphics} % for pdf, bitmapped graphics files
\usepackage{epsfig} % for postscript graphics files
\usepackage{mathptmx} % assumes new font selection scheme installed
\usepackage{times} % assumes new font selection scheme installed
\usepackage{amsmath} % assumes amsmath package installed
\usepackage{amssymb}  % assumes amsmath package installed
\usepackage{makeidx}
\usepackage{float}
\usepackage{balance}
\usepackage{booktabs}
\usepackage{bbm}
\usepackage{mathrsfs}
\usepackage{enumerate}
\usepackage{multicol}
\usepackage{algorithm}
\usepackage{algpseudocode}
\usepackage[nolist]{acronym} %acronyms
\usepackage{mathrsfs}  
\usepackage[normalem]{ulem}
\usepackage{balance}
\usepackage{xr}
\usepackage{xcolor}
\usepackage{dsfont}
\usepackage{amsthm}

\usepackage{cleveref}

\hypersetup{
  colorlinks = true,
  allcolors = siaminlinkcolor,
  urlcolor = siamexlinkcolor,
}
\colorlet{siaminlinkcolor}{green!50!black}
\colorlet{siamexlinkcolor}{red!50!black}

\usepackage{environ} %for making custom environments

\usepackage{graphicx}
\usepackage{textcomp}

\newcommand{\beginsupplement}{
        \setcounter{table}{0}
        \renewcommand{\thetable}{S\arabic{table}}
        \setcounter{figure}{0}
        \renewcommand{\thefigure}{S\arabic{figure}}
        \setcounter{section}{0}
        \renewcommand{\thesection}{S\arabic{section}}
        \setcounter{equation}{0}
        \renewcommand{\theequation}{S\arabic{equation}}
        \setcounter{algorithm}{0}
        \renewcommand{\thealgorithm}{S\arabic{algorithm}}
     }

\newcommand{\rd}{{\mathrm d}}

\newcommand{\rvec}{{\mathrm{vec}}}

\newcommand{\vx}{{\bf x}}

\newcommand{\calL}{{\cal L}}

\newcommand{\calU}{{\cal U}}

\newcommand{\Rb}{\mathbb{R}}

\newcommand{\Nb}{\mathbb{N}}

\newcommand{\Udbar}{\bar{U}_d}
\newcommand{\Ubbar}{\bar{U}_b}
\newcommand{\Xbar}{\bar{X}}

\newcommand{\delUb}{\delta U_b}
\newcommand{\delUd}{\delta U_d}
\newcommand{\delU}{\delta U}
\newcommand{\delX}{\delta X}

\newcommand{\ddt}{\frac{\rd}{\rd t}}
\newcommand{\half}{\frac{1}{2}}

\newcommand{\dSxi}{\rd S_\xi}
\newcommand{\dSeta}{\rd S_\eta}

\NewEnviron{eqsp}[1][]{%
\begin{equation} \label{#1}
\begin{split}
\BODY
\end{split}
\end{equation}}
\newcommand{\vm}{{\bf m}}

\newcommand{\Langle}{{\Big\langle}}
\newcommand{\Rangle}{{\Big\rangle}}

\newcommand{\limiinfty}{{\lim_{i\rightarrow\infty}}}

\makeatletter
\renewcommand{\maketag@@@}[1]{\hbox{\m@th\normalsize\normalfont#1}}%
\makeatother

\newtheorem{theorem}{Theorem}[section]
\newtheorem{corollary}{Corollary}[section]
\newtheorem{lemma}{Lemma}[section]
\newtheorem{prop}{Proposition}[section]
\newtheorem{remark}{Remark}

\newtheorem{assumption}{Assumption}

\catcode`\_=11\relax
    \newcommand\email[1]{\_email #1\q_nil}
    \def\_email#1@#2\q_nil{%
      \href{mailto:#1@#2}{{\emailfont #1\emailampersat #2}}
    }
    \newcommand\emailfont{}
    \newcommand\emailampersat{\small@}
    \catcode`\_=8\relax

\begin{acronym}
\acro{DNN}{Deep Neural Network}
\acro{ODE}{Ordinary Differential Equation}
\acro{SPDE}{Stochastic Partial Differential Equation}
\acro{FNN}{Feed-forward Neural Network}
\acro{CNN}{Convolutional Neural Network}
\acro{DP}{Dynamic Programming}
\acro{LSTM}{Long-Short Term Memory}
\acro{FC}{Fully Connected}
\acro{DDP}{Differential Dynamic Programming}
\acro{HJB}{Hamilton-Jacobi-Bellman}
\acro{PDE}{Partial Differential Equation}
\acro{LQR}{Linear Quadratic Regulator}
\acro{RDE}{Riccati Differential Equation}
\acro{PI}{Path Integral}
\acro{NN}{Neural Network}
\acro{SOC}{Stochastic Optimal Control}
\acro{RL}{Reinforcement Learning}
\acro{MPC}{Model Predictive Control}
\acro{IL}{Imitation Learning}
\acro{RNN}{Recurrent Neural Network}
\acro{DL}{Deep Learning}
\acro{RN}{Radon-Nikodym}
\acro{SGD}{Stochastic Gradient Descent}
\acro{SDE}{Stochastic Differential Equation}
\acro{VRL}{Variational Reinforcement Learning}
\acro{IDVRL}{Infinite Dimensional Variational Reinforcement Learning}
\acro{1D}{1-dimensional}
\acro{2D}{2-dimensional}
\acro{3D}{3-dimensional}
\acro{ROM}{Reduced Order Model}
\acro{STSO}{Spatio-Temporal Stochastic Optimization}
\acro{ANN}{Artificial Neural Network}
\acro{ADPL}{Actuator Design and Policy Learning}
\acro{STDDP}{Spatio-Temporal \ac{DDP}}
\end{acronym}

\begin{document}
	
	\title{
	\rule{\linewidth}{4pt}\vspace{0.3cm} \Large \textbf{
	  Spatio-Temporal Differential Dynamic Programming \\
	  for Control of Fields 
	}\\ \rule{\linewidth}{1.5pt}}
	\author{Ethan N. Evans$^{a,}\thanks{Corresponding Author. Email: \email{eevans89@gmail.com}}~$, Oswin So$^{b}$, Andrew P. Kendall$^{a}$, Guan-Horng Liu$^a$,\\ and Evangelos A. Theodorou$^{a,c}$\\ \vspace{-0.1cm}
	\small{$^a$Georgia Institute of Technology, Department of Aerospace Engineering} \\ \vspace{-0.2cm}
	\small{$^b$Georgia Institute of Technology, College of Computing} \\ \vspace{-0.2cm}
	\small{$^c$Georgia Institute of Technology, Institute of Robotics and Intelligent Machines} }
	
	\date{\small{This manuscript was compiled on \today}}
	
	\maketitle

\begin{abstract}
We consider the optimal control problem of a general nonlinear spatio-temporal system described by Partial Differential Equations (PDEs). Theory and algorithms for control of spatio-temporal systems are of rising interest among the automatic control community and exhibit numerous challenging characteristic from a control standpoint. Recent methods focus on finite-dimensional optimization techniques of a discretized finite dimensional ODE approximation of the infinite dimensional PDE system. In this paper, we derive a differential dynamic programming (DDP) framework for distributed and boundary control of spatio-temporal systems in infinite dimensions that is shown to generalize both the spatio-temporal LQR solution, and modern finite dimensional DDP frameworks. We analyze the convergence behavior and provide a proof of global convergence for the resulting system of continuous-time forward-backward equations. We explore and develop numerical approaches to handle sensitivities that arise during implementation, and apply the resulting STDDP algorithm to a linear and nonlinear spatio-temporal PDE system. Our framework is derived in infinite dimensional Hilbert spaces, and represents a discretization-agnostic framework for control of nonlinear spatio-temporal PDE systems.

\end{abstract}

%===============================================================================

\section{Introduction}\label{sec:intro}
Many complex natural processes are governed by systems of equations with spatio-temporal dependence, and are typically described by \ac{PDE}. These systems are ubiquitous in nature and can be found in most disciplines of engineering and applied physics. The range of natural processes includes fluid flow governed by the Navier-Stokes equation, sub-atomic particle systems governed by the Schrodinger equation, activation of neurons governed by the Nagumo equation~\cite{lord_powell_shardlow_2014}, and flame front propagation in combustion systems governed by the Kuramoto-Sivashinsky equation~\cite{gomes2017controlling}.

Despite their ubiquity in nature and engineering, theory and numerical methods for control of spatio-temporal systems remains challenging due to the time-delay, dramatic under-actuation, high system dimensionality, and multi-modal bifurcations, which are often inherent in their dynamics. Furthermore, existence and uniqueness of solutions remains an open problem for many systems, and when they do exist, they typically only have a weak notion of differentiability. Analysis of their performance must be treated with calculus over functionals, and their state vectors are often described by vectors in an infinite-dimensional time-indexed Hilbert space even for scalar \ac{1D} \acp{PDE}. Put together, mathematically consistent and numerically realizable algorithms for control of spatio-temporal systems represents many of the largest current-day challenges facing the automatic control community.

The majority of recent methods for control of spatio-temporal systems typically reduce \acp{PDE} into a finite set of \acp{ODE} through \acp{ROM}, and apply standard finite-dimensional optimization methods which result in algorithms specific to the \ac{ROM} used. Within this paradigm, deep learning methods have successfully been applied on policy networks in the finite dimensional setting for controlling Navier-Stokes systems \cite{rabault2019artificial,bieker2019deep,mohan2018deep,nair2019cluster}, for soft robotic systems \cite{satheeshbabu2019open,spielberg2019learning}, as well as for many other systems \cite{farahmand2017deep}. These methods are often specific to a discretization scheme and represent a \textit{discretize-then-optimize} approach. Some such methods can introduce new phenomena in the latent space represetnation, as in \cite{morton2018deep}, where the resulting deep Koopman approach can be shown to violate linear stabilizability conditions of the latent space dynamics.

External to the machine-learning literature are infinite-dimensional methods found in the control theory literature~\cite{lasiecka2000control,troltzsch2010optimal}, which are dominated by linear or linearization-based approaches, which include \ac{LQR} approaches for linear \acp{PDE}, and forward-backward approaches, which include approaches due to the Pontryagin Maximum Principle (PMP) \cite{troltzsch2010optimal,sumin2009first,yong1992pontryagin}. Indeed local linearization methods allow for optimal solutions of an approximate problem, however require knowledge of linearization points a-priori. On the other hand, forward-backward schemes provide a nominal trajectory and optimization-based control update scheme at the expense of the backpropagation of a coupled system equation.

In contrast to Pontryagin methods which yield a state-independent backward equation and an open-loop controller, methods founded on the Bellman principle of optimality utilize backward equations that are state-dependent and yield closed-loop control solutions. Methods such as \ac{DDP} have decades of established history in the finite dimensional automatic control literature. Modern variations include control limits \cite{tassa2014control}, state constraints \cite{aoyama2020constrained}, receding horizons \cite{tassa2007receding}, belief space control \cite{pan2014probabilistic,pan2018efficient}, game-theoretic control \cite{sun2015game}, control on Lie groups \cite{boutselis2019numerical}, and using polynomial chaos variational integrators \cite{boutselis2016stochastic}.

A previous attempt exists to extend the \ac{DDP} framework to spato-temporal systems in infinite dimensions \cite{tzafestas1969differential}, however this approach has several flaws and mathematical inconsistencies, as pointed out in \cite{sakawa1972matrix}. Additionally, the \ac{DDP} method has had significant growth since the early works \cite{jacobson1970}. Decades of advancement include linearization around the nominal trajectory as opposed to the optimal trajectory which decreases sensitivities of convergence behavior to the initial conditions, regularization in the second order backward equation to increase numerical stability, treatment of state and control constraints, and optimization over time horizon.

In light of the apparent literature gap, this manuscript is devoted to the development of \ac{DDP} methods for spatio-temporal systems in infinite dimensions. Specifically, we derive the \ac{STDDP} framework incorporating modern theoretical techniques, we demonstrate that the resulting system of forward-backward equations generalizes both the \ac{LQR} solution in infinite dimensions and \ac{DDP} in finite dimensions, we provide a proof of convergence for the resulting system of continuous-time forward-backward equations, we explore and develop numerical approaches to handle sensitivities that arise due to discretization, and apply the resulting algorithm to linear and nonlinear spatio-temporal \ac{PDE} systems. In contrast to recent machine learning methods, our optimization is developed entirely in Hilbert spaces, and represents an \textit{optimize-then-discretize} approach. As a result, the framework is a continuous-time formulation which is \textit{agnostic} to discretization scheme during implementation.

%===============================================================================

\section{Preliminaries and Problem Statement}

Let $D \subseteq \Rb^n$ denote a measurable connected open domain of $\Rb^n$ describing the space on which the system evolves. Let $S \subseteq \Rb^n$ denote the boundary of $D$, let $\bar{D}$ denote the closure of the domain, i.e. $\bar{D} = D \cup S$, and let $T = [t_0, t_f]$ denote some arbitrary time domain. In fields representation, a general form of a deterministic \ac{PDE} dynamical system is given by
\begin{align}
    \partial_t X(t,x) &= F(t,x,X(t,x), U_d(t,x)), \quad  x \in D \label{eq:field_dynamics} \\
    0 &= N(t,x,X(t,x),U_b(t,x)), \quad x \in S \label{eq:field_boundary}\\
    X(t_0,x) &= X_0(x), \quad x \in \bar{D} \label{eq:field_IC},
\end{align}
where $X:T \times D \rightarrow \Rb^n$ is the state. This problem has two measurable control functions, $U_b: T \times S \rightarrow \Rb^l$ which correspond to actuation on the boundary, and $U_d: T \times D \rightarrow \Rb^k$ which corresponds to actuation distributed throughout the field excluding the boundary. The dynamics evolve by some measurable functional $F: T \times D 
\times \Rb^n \times \Rb^k \rightarrow \Rb^n$ that is potentially nonlinear in the state function $X(t,x)$ or the control function $U_d(t,x)$, with a boundary condition functional $N: T \times S \times \Rb^n \times \Rb^l \rightarrow \Rb^n$ that is also potentially a nonlinear functional of the state or control functions, and can be any type of boundary condition (e.g. Neumann, Dirichlet, etc.). 

We can equivalently write \cref{eq:field_dynamics,eq:field_boundary,eq:field_IC} in the time-indexed Hilbert spaces perspective by first properly defining Hilbert spaces, as in \cite{sakawa1972matrix}. Let $L_2^n(D)$ denote the Hilbert space of $n$-vector functions square integrable over $D$ with inner product 
\begin{equation}
    \Langle X_1, X_2 \Rangle = \int_D X_1^\top(x) X_2(x) \rd x,
\end{equation}
where $\rd x = \rd x_1 \rd x_2 \cdots \rd x_n$ is shorthand notation for the generalized volume integration over $\Rb^n$. This is the Hilbert space of the domain, and we similarly define the Hilbert space over the boundary. Let $L_2^n(S)$ denote the Hilbert space of $n$-vector functions square integrable over $S$ with inner product
\begin{equation}
    \Langle X_1, X_2 \Rangle_S = \int_S X_1(\xi) X_2(\xi) \dSxi,
\end{equation}
where $\dSxi$ is an infinitesimal surface element of the boundary at a point $\xi \in S$. Let $\calL(U,V)$ denote the space of linear bounded operators from $U$ into $V$. If we regard $X(t,x)$ as an element of $L_2^n(D)$, then we can rewrite \cref{eq:field_dynamics,eq:field_boundary,eq:field_IC} as
\begin{align}
    \ddt X(t) &= F\big(t,X(t),U_d(t)\big), \quad X \in L_2^n(D), \;\; t\in T \label{eq:hilbert_dynamics} \\ 
    0 &= N\big(t, X(t),U_b(t) \big), \quad X \in L_2^n(S), \;\; t\in T \label{eq:hilbert_boundary}\\
    X(t_0) &= X_0, \label{eq:hilbert_IC}
\end{align}
where $X(t), X_0 \in L_2^n(D)$ are respectively the Hilbert space state vector and initial conditions, $U_d(t) \in L_2^k(D)$ is the Hilbert space distributed control vector, $U_b \in L_2^l(S)$, is the Hilbert space boundary control vector, $F : T \times  L_2^n(D) \times L_2^k(D) \rightarrow L_2^n(D)$ is a potentially nonlinear measurable function on the domain Hilbert space, and $N: T \times L_2^n(D) \times L_2^l(S) \rightarrow L_2^n(S)$ is a potentially nonlinear measurable function on the boundary Hilbert space.

\begin{remark}
The field functional perspective of \cref{eq:field_dynamics,eq:field_boundary,eq:field_IC} and the time-indexed Hilbert space perspective of \cref{eq:hilbert_dynamics,eq:hilbert_boundary,eq:hilbert_IC} are consistent in the sense that they share identical solutions up to the transformation between the perspectives used above.
\end{remark}

\begin{assumption}\label{assumption_existence}
    The \ac{PDE} system in fields representation given by \cref{eq:field_dynamics,eq:field_boundary,eq:field_IC} is well posed in the sense of Hadamard, and admits a unique weak solution $X(t,x)$, $t \in T, x \in \bar{D}$ for each initial condition $X_0(x) \in \Rb^n$.
\end{assumption}
Depending on the specific form of the \ac{PDE}, this assumption can have varying degrees of severity, however in general it is a mild assumption. Please refer to \cite{evans1997partial} for more details on existence and uniqueness of various \acp{PDE}. Despite the potential severity, it is an assumption that is required henceforth. 
% Existence and uniqueness of solutions of the Hilbert space form are given by the following remark.
\begin{remark}
    If \cref{assumption_existence} holds, then the \ac{PDE} system in Hilbert space representation given by \cref{eq:hilbert_dynamics,eq:hilbert_boundary,eq:hilbert_IC} is also well posed in the sense of Hadamard, and admits a unique weak Hilbert space solution $X(t) \in L_2^n(\bar{D})$, $t \in T$ for each Hilbert space initial condition $X_0 \in L_2^n(\bar{D})$.
\end{remark}

This remark has an obvious proof (e.g. by contradiction) that is omitted. Throughout this work, we go back and forth between these two notational perspectives: the spatially varying fields perspective, and the time-indexed Hilbert space perspective. While the fields perspective demonstrates the spatial integration that is central to the Volterra-Taylor expansions more clearly, the time-indexed Hilbert space perspective will often yield a more compact notation that is easier to treat with familiar algebraic operations. Whenever we suppress the dependencies on the spatial variable $x$, the variables are assumed to be in time-indexed Hilbert spaces. 

In order to arrive at the optimal control problem, we first define the measurable cost functional in fields representation as
\begin{equation}\label{eq:fields_cost}
\begin{split}
    J\big(t, X(t,x), U_d(t,x), U_b(t,x) \big) := \phi\big(t_f, X(t_f,x) \big) + \int_{t_0}^{t_f} L\big(t, X(t,x), U_d(t,x), U_b(t,x) \big) \rd t,
\end{split}
\end{equation}
where $\phi: T \times \Rb^n \rightarrow \Rb$ is some measurable real-valued terminal cost functional, and $L: T \times \Rb^n \times \Rb^k \times \Rb^l \rightarrow \Rb$ is a measurable real-valued running cost functional. In time-indexed Hilbert spaces, the cost functional becomes
\begin{equation}\label{eq:hilbert_cost}
\begin{split}
    J\big(t, X(t), U_d(t), U_b(t) \big) = \phi\big(t_f, X(t_f) \big)  + \int_{t_0}^{t_f} L\big(t, X(t), U_d(t), U_b(t) \big) \rd t,
\end{split}
\end{equation}
where $J:T \times L_2^n(D) \times L_2^k(D) \times L_2^l(S) \rightarrow \Rb$, $\phi: T \times L_2^n(D) \rightarrow \Rb$, and $L: T \times L_2^n(D) \times L_2^k(D) \times L_2^l(S) \rightarrow \Rb$ are the equivalent measurable real-valued functionals in Hilbert spaces. The value functional is defined as
\begin{align}\label{eq:val_fcn}
    V\big(X(t), t \big) := \min_{U_d, U_b} \Big[ J\big( X(t), U_d(t), U_b(t) \big) \Big].
\end{align}
Due to the Bellman Principle of Optimality, one can form the \ac{HJB} equation as \cite{tzafestas1969differential,sakawa1972matrix}
\begin{align}\label{eq:HJB}
    -\partial_t V\big(X(t),t\big) &= \min_{U_d, U_b} \Big[  L(t,X,U_d,U_b) + \Langle V_X(t,X) , F(t, X, U_d) \Rangle  \Big], \\
    V\big(X(t_f), t_f\big) &= \phi\big(t_f, X(t_f)\big) =: V_f \in \Rb,
\end{align}
where we write $\partial_t = \frac{\partial}{\partial t}$ to denote the normal partial derivative of a function with respect to a variable, and use subscript $X$, $U_b$, or $U_d$ to denote the Gateaux partial derivative of a functional or operator with respect to an operator function. One can carry out the same derivation using Volterra's notion of functional derivative \cite{volterra1959theory}. Note that $V(X(t),t)$ is a function of time, and a functional of $X(t)$. Also, it should be noted that the \ac{HJB} equation in \cref{eq:HJB} is a backwards nonlinear \ac{PDE}.

\begin{assumption}\label{assumption_HJB_existence}
    The backwards \ac{PDE} in \cref{eq:HJB} admits a unique viscosity solution $V\big(t,X(t)\big)$, $t \in T$, $X(t) \in L_2^n(\bar{D})$ for each terminal condition $V\big(X(t_f), t_f\big) = V_f := \phi\big(t_f, X(t_f)\big) \in \Rb$.
\end{assumption}

The \ac{DDP} framework solves the \ac{HJB} equation in \cref{eq:HJB} iteratively via expansions of the value functional, cost functional, dynamics operator function, and boundary operator function to given order. Typically, the value functional and cost functional are expanded to second order so that the resultig \ac{HJB} becomes a quadratic optimization problem with a unique optimal control minimizer.

Quadratic expansions also allow for proofs of global convergence and even proofs of quadratic convergence, that in finite dimensions, initially relied on well known convergence properties of the Newton method of optimization \cite{shoemaker1990proof,yakowitz1984computational} for quadratic problems. Under similar reasoning, the dynamics are typically either expanded to first or second order.

%===============================================================================

\section{Expansions of the Cost, Value, Field, and Boundary}

The approach in this paper is a spatio-temporal \ac{DDP} approach that is analogous to the finite dimensional \ac{DDP} apparoch of \cite{sun2014continuous}. Therein, the authors discuss the fundamental differences between their derivation, and the original derivation by Jacobson and Mayne \cite{jacobson1970}. The derivation by Jacobson and Mayne, of which a similar flavor is followed in \cite{tzafestas1969differential}, is based on the restrictive assumption that the nominal control trajectory $\bar{u}$ is sufficiently close to the optimal control solution $u^*$. This is circumvented by performing expansions around a nominal trajectory. Define a nominal state and control triple $(\Xbar, \Udbar, \Ubbar)$ and the variations $\delX := X - \Xbar$, $\delUd := U_d - \Udbar$, and $\delUb := U_b - \Ubbar$. In order to properly write the expansions, we require the following assumption

\begin{assumption}\label{assumption_differentiable}
    The dynamics function $F$ and boundary function $N$ are differentiable almost everywhere, the running cost functional $L$ and terminal cost functional $\phi$ are twice differentiable almost everywhere, and the value functional $V$ is three times differentiable almost everywhere. 
    These stated derivatives are defined in the Gateaux sense with respect to the state and control triple $(X, U_d, U_b)$, and are square integrable in the Lebesgue sense. That is, the stated Gateaux derivative of each functional exists $\forall (X, U_d, U_b) \in \big(L_2^n(D), L_2^k(D), L_2^l(S)\big)$ except on a properly defined set of measure zero.
\end{assumption}

As previously stated, the value functional is a function of time $t$, but a functional of the spacetime function $X(t,x)$. Thus the value functional is expanded via a Volterra-Taylor functional expansion \cite{volterra1959theory}
\begin{eqsp}[eq:V_expansion_fields]
    V\big(t,\Xbar(t,x) + \delX(t,x)\big) &= V\big(t, \Xbar(t,x)\big) + \int_D V_X^\top\big(t,\Xbar(t,x)\big) \delX(t,x) \rd x \\
    &\quad + \half \int_D \int_D \delX^\top(t,x)  V_{XX}(t,x,y) \delX(t,y) \rd x \rd y + O(\delta^3).
\end{eqsp}

We maintain connection to the Hilbert space perspective by defining Hilbert space operators for each kernel function. Define the operator $V_{XX}(t,X) \in \calL \big(L_2^n(D), L_2^n(D)\big)$ as 
\begin{equation}
    V_{XX}\big(t, X\big)W(t) := \int_D V_{XX}(t,x,y) W(t,y) \rd y,
\end{equation}
where $V_{XX}(t,x,y)$ is the kernel function. In order to form the left-hand side of the \ac{HJB} \cref{eq:HJB}, we apply a re-arranged definition of the total differential \cite{volterra1959theory}, given by
\begin{align}\label{eq:total_deriv}
    \partial_t (\cdot) = \frac{\rd}{\rd t} (\cdot) - \Langle  (\cdot)_X , F(t, X, U_d) \Rangle,
\end{align}
which holds for any functional that explicitly depends on $X$ and $t$. In order to simplify notation, we suppress arguments when  functionals are evaluated on the nominal trajectory triple. We apply \cref{eq:total_deriv} to each term on the right-hand side of \cref{eq:V_expansion_fields}, to yield the left-hand side of the \ac{HJB}, which in Hilbert spaces, has the form
\begin{align}\label{eq:HJB_LHS}
    -\partial_t V(t,\Xbar + \delX) &= -\ddt \bigg(V +  \Langle  V_X, \delX \Rangle  + \half \Langle \delX , V_{XX} \delX \Rangle  \bigg) + \Langle V_X, F \Rangle + \Langle V_{XX} F, \delX \Rangle  \\
    &\quad + \half  \Langle F V_{XXX} \delX', \delX \Rangle, 
\end{align}
where for the third order Gateaux derivative $V_{XXX}$, we have defined the tensor operator in time-indexed Hilbert spaces $V_{XXX}(t,X) \in \calL\big(L_2^n(D) \times L_2^n(D), L_2^n(D)\big)$ as
\begin{equation}
    U(t) V_{XXX}\big(t, X(t)\big) W(t) := \int_D \int_D U^\top(t,x) V_{XXX}(t,x,y,z) W(t,y) \rd x \rd y.
\end{equation}
The 4-D kernel function $V_{XXX}(t,x,y,z)$ is assumed to be symmetric about all three spatial axes for simplicity.

Next, we expand the cost functional with a Volterra-Taylor expansion to second order, which in time-indexed Hilbert spaces has the form
\begin{eqsp}[eq:L_expand]
    L(t,\Xbar + \delX, \Udbar + \delUd, \Ubbar + \delUb) 
    &= L + \Langle L_X,\, \delX \Rangle + \Langle L_{U_d}, \delUd \Rangle + \Langle L_{U_b},\, \delUb \Rangle_S + \half \Langle \delX, L_{XX} \delX' \Rangle   \\
    &\quad + \half \Langle \delUd, \Big(L_{U_d X} + L_{X U_d}^\top \Big) \delX' \Rangle  + \half \Langle \delUb, \Big(L_{U_b X} + L_{X U_b}^\top \Big) \delX' \Rangle_S  \\
    &\quad + \half\Langle \delUd, L_{U_d U_d} \delUd' \Rangle + \half \Langle \delUb, L_{U_b U_b} \delUb' \Rangle_S + O(\delta^3),
\end{eqsp}
where we have defined the operators
\begin{align*}
    L_{XX}\big(t, X(t),U_b(t),U_d(t)\big) W(t) &:= \int_D L_{XX} (t,x,y) W(t, y) \rd y \\
    L_{X U_d}\big(t, X(t),U_b(t),U_d(t)\big) W(t) &:= \int_D L_{X U_d} (t,x,y) W(t, y) \rd y\\
    L_{X U_b}\big(t, X(t),U_b(t),U_d(t)\big) W(t) &:= \int_S L_{X U_b} (t,\xi,\eta) W(t,\eta) \rd S_\eta\\
    % L_{U_d U_b}\big(t, X(t),U_b(t),U_d(t)\big) A(t) &:= \int_S L_{U_d U_b} \big(X(\xi,t),X(\eta,t), t\big) A(\eta,t) \rd S_\eta \textcolor{red}{\text{ verify this one}} \\
    L_{U_d U_d}\big(t, X(t),U_b(t),U_d(t)\big) W(t) &:= \int_D L_{U_d U_d} (t,x,y) W(t, y) \rd y \\
    L_{U_b U_b}\big(t, X(t),U_b(t),U_d(t)\big) W(t) &:= \int_S L_{U_b U_b}(t,\xi,\eta) W(t,\eta) \rd S_\eta,
\end{align*}
and similarly defined operators for $L_{U_d X}$, $L_{U_b X}$. 

\begin{assumption}
The measurable kernel functions $L_{XX}$, $L_{X U_d}$, $L_{X U_b}$ are spatially symmetric and positive semi-definite. The measurable kernel functions $V_{XX}$, $L_{U_d U_d}$, $L_{U_b U_b}$ are spatially symmetric and positive definite. The omitted cross term operators $L_{U_b U_d}$ and $L_{U_d U_b}$ are null operators.
\end{assumption}

Note the assumption that cross terms between boundary and distributed control (i.e. $L_{U_b U_d}$ and $L_{U_d U_b}$) are zero. This is a fairly benign assumption since cost functionals are often composed of pure quadratics in either $U_d$ or $U_b$, but not both. Including these cross terms also yields optimal update equations for boundary and distributed control that are coupled to each other, and thus impose mathematical and implementation difficulties.

Next, the dynamics and boundary are expanded around the nominal trajectory. The dynamics functional $F\big(t, X(t), U_d(t)\big)$ and boundary functional $N\big(t, X(t), U_b(t) \big)$ map into $L_2^n(D)$ and $L_2^n(S)$, respectively, and are not real-valued functionals, so it is appropriate to treat them as operator functions instead of as functionals despite having explicit dependence on functions $\Xbar, \Udbar, \Ubbar$. In Hilbert space notation, the operator Taylor expansion of the dynamics and boundary have the form
\begin{align} 
    F(t, \Xbar + \delX,\Udbar + \delUd) &= F(t, \Xbar, \Udbar) + F_X^\top(t, \Xbar, \Udbar) \delX + F_{U_d}^\top(t, \Xbar, \Udbar) \delUd + O(\delta^2) \label{eq:F_expand} \\
    N(t, \Xbar + \delX, \Ubbar + \delUb) &= N(t, \Xbar, \Ubbar) + N_X^\top(t, \Xbar, \Ubbar) \delX  + N_{U_b}^\top(t, \Xbar, \Ubbar) \delUb + O(\delta^2) \label{eq:N_expand}
\end{align}
We obtain the right-hand side of the \ac{HJB} \cref{eq:HJB} by plugging \cref{eq:L_expand,eq:F_expand,eq:N_expand}, and a Volterra-Taylor expansion of $V_X$. After simplification, the right-hand side of the \ac{HJB} \cref{eq:HJB} becomes
\begin{eqsp}[eq:HJB_RHS_simp]
    \min_{\delUd, \delUb} \bigg[& L + \Langle L_X,\, \delX \Rangle + \Langle L_{U_d}, \delUd \Rangle + \Langle L_{U_b},\, \delUb \Rangle_S + \half \Langle \delX,  L_{XX}\, \delX' \Rangle \\
    & + \half \Langle \delUd,\, \Big( L_{U_d X} + L_{X U_d}^\top \Big) \delX \Rangle + \half \Langle \delUb,\, \Big(L_{U_b X} + L_{X U_b}^\top \Big) \delX' \Rangle_S  + \half \Langle \delUd, L_{U_d U_d} \,\delUd' \Rangle \\
    & + \half \Langle \delUb,\, L_{U_b U_b} \,\delUb' \Rangle_S + \Langle V_X,\, F \Rangle + \Langle V_X,\, F_X^\top \delX \Rangle + \Langle V_X,\, F_{U_d}^\top \delUd \Rangle + \Langle \delX,\, V_{XX} F  \Rangle \\
    &+ \Langle \delX,\, V_{XX} F_X^\top \delX' \Rangle + \Langle \delX,\, V_{XX} F_{U_d}^\top \delUd \Rangle + \half \Langle \delX,\, V_{XXX}F \delX' \Rangle \bigg].
\end{eqsp}
Equating \cref{eq:HJB_LHS} to \cref{eq:HJB_RHS_simp} and canceling common terms yields 
\begin{eqsp}[eq:HJB_LHS_RHS]
    -\ddt \bigg(V +   \Langle  V_X, &\delX \Rangle  + \half \Langle \delX , V_{XX} \delX' \Rangle \bigg) \\
    = \min_{\delUd, \delUb} \bigg[& L + \Langle L_X,\, \delX \Rangle + \Langle L_{U_d}, \delUd \Rangle + \Langle L_{U_b},\, \delUb \Rangle_S  + \half \Langle \delX,  L_{XX}\, \delX' \Rangle \\
    &+ \half \Langle \delUd,\, \Big( L_{U_d X} + L_{X U_d}^\top \Big) \delX \Rangle + \half \Langle \delUb,\, \Big(L_{U_b X} + L_{X U_b}^\top \Big) \delX' \Rangle_S + \half \Langle \delUd, L_{U_d U_d} \,\delUd' \Rangle \\
    & + \half \Langle \delUb,\, L_{U_b U_b} \,\delUb' \Rangle_S + \Langle V_X,\, F_X^\top \delX \Rangle + \Langle V_X,\, F_{U_d}^\top \delUd \Rangle + \Langle \delX,\, V_{XX} F_X^\top \delX' \Rangle \\
    & + \Langle \delX,\, V_{XX} F_{U_d}^\top \delUd \Rangle \bigg].
\end{eqsp}
% The intermediate steps to obtain \cref{eq:HJB_RHS_simp} and \cref{eq:HJB_LHS_RHS} can be found in the Supplemental Material \textbf{specify section}.
%\guan{should $\Langle \delX,\, V_{XX} F_X^\top \delX' \Rangle$ in \cref{eq:HJB_LHS_RHS} be $\Langle  V_{XX} \delX,\,F_X^\top \delX' \Rangle$ if you are substituting \cref{eq:green_HJB_term2} later to get \cref{eq:HJB_LHS_RHS_green} }

\begin{remark}
The exact singleton Newton minimizer $\delUb^*$ of the approximate \ac{HJB} equation \cref{eq:HJB_LHS_RHS} does not incorporate the value functional $V(t,\Xbar)$ or its derivatives $V_X(t,\Xbar)$, $V_{XX}(t,\Xbar)$.
% The approximate HJB equation \cref{eq:HJB_LHS_RHS} is not quadratic in $\delUb$. Futhermore, \cref{eq:HJB_LHS_RHS} is not strictly convex in $\delUb$. Thus, \cref{eq:HJB_LHS_RHS} does not have an exact singleton Newton optimization solution for $\delUb$. %Thus Newton optimization cannot be performed exactly to find a singleton minimizer.
\end{remark}

This is an important point. The value functional is defined as the minimization surface of the original problem in \cref{eq:val_fcn} and the apparent decoupling between the optimal update $\delUb^*$ and the value functional and/or its derivatives within the resulting approximate \ac{HJB} \cref{eq:HJB_LHS_RHS} yields a naive update. The authors in \cite{tzafestas1969differential} and \cite{tzafestas1968optimal} realize this fact, and use the Green's theorem in order to incorporate boundary information into specific terms in \cref{eq:HJB_LHS_RHS}. However, there are errors in their application of Green's theorem in the multivariate case, as noted in \cite{sakawa1972matrix}.

%===============================================================================

\section{Green's Theorem in Hilbert Spaces}

Green's theorem is used widely in calculus to relate the volume integral of the interior of a region to a surface integral of its boundary. In the context of \ac{STDDP}, it allows us to capture pertinent effects of the value function on the boundary.

\vspace{0.5em}
\begin{assumption}\label{assumption_greens}
    $F_X$ is a linear operator with standard form given by $A_x(t,x)$ in \cite{sakawa1972matrix}, and $N_X$ is a linear operator with standard form given by $\beta_A(t,\xi)$ in \cite{sakawa1972matrix}.
\end{assumption}

% With this, we provide a generalized version of Green's theorem in Hilbert spaces for clarity. 
% More details can be found in Supplemental Material \textbf{specify section}.

\vspace{0.5em}
\begin{theorem}\label{thm:green}
Let $Y(t), Z(t) \in L_2^n(\bar{D})$. Under \cref{assumption_greens}, the following holds:
\begin{equation}
\begin{split}
    &\Langle Y(t), F_X\big(t,X(t), \Udbar(t) \big) Z(t) \Rangle - \Langle Z(t), F_X^*\big(t,X(t), \Udbar(t) \big) Y(t) \Rangle \\
    &= \Langle Y(t), N_X\big(t,X(t), \Ubbar(t)\big) Z(t) \Rangle_S - \Langle Z(t), N_X^*\big(t,X(t), \Ubbar(t) Y(t)\big)\Rangle_S
\end{split}
\end{equation}
\end{theorem}
\vspace{1em}

The equivalent fields representation can be found in \cite{sakawa1972matrix}, and the proof is a standard result (c.f. \cite{friedman2008partial}). The following corollary is a direct application of \cref{thm:green} to the applicable terms of the \ac{HJB} in \cref{eq:HJB_LHS_RHS}.

\begin{corollary}
    If \cref{assumption_greens} holds, then
    \begin{align}
         \Langle V_X, F_X^\top \delX \Rangle &= \Langle  \delX, F_X^* V_X \Rangle  - \Langle V_X, \Delta N \Rangle_S -  \Langle V_X, N_{U_b}^\top \delUb \Rangle_S - \Langle \delX, N_X^* V_X \Rangle_S\;\;, \label{eq:green_HJB_term1}
    \end{align}
    and 
    % \guan{for the RHS of \cref{eq:green_HJB_term2}, should $\delX'$ and $\delX$ switch with each other according to Thm 4.1?}
    \begin{align}
        \Langle  V_{XX} \delX,\,F_X^\top \delX' \Rangle
        &= \Langle \delX ,\,F_X^*  V_{XX} \delX' \Rangle - \Langle V_{XX} \delX,\, \Delta N \Rangle_S   - \Langle V_{XX} \delX, N_{U_b}^\top \delUb \Rangle_S  - \Langle \delX,\, N_X^* V_{XX} \delX' \Rangle_S \;\;, \label{eq:green_HJB_term2}
\end{align}
    where $\Delta N = N(X + \delX, U_b + \delUb) - N(X,U_b)$, $F_X^*$ is the adjoint operator of $F_X$, $N_X^*$ is the adjoint operator of $N_X$, and we have suppressed explicit time dependencies for simplicity.
\end{corollary}

Plugging equations \cref{eq:green_HJB_term1,eq:green_HJB_term2} into \cref{eq:HJB_LHS_RHS} yields
\begin{eqsp}[eq:HJB_LHS_RHS_green]
    -\ddt \bigg(V + \Langle  V_X, &\delX \Rangle  + \half \Langle \delX , V_{XX} \delX' \Rangle \bigg) \\
    = \min_{\delUd, \delUb} \bigg[& L + \Langle L_X,\, \delX \Rangle + \Langle L_{U_d}, \delUd \Rangle + \Langle L_{U_b},\, \delUb \Rangle_S  + \half \Langle \delX,  L_{XX}\, \delX' \Rangle \\
    & + \half \Langle \delUd,\, \Big( L_{U_d X} + L_{X U_d}^\top \Big) \delX' \Rangle + \half \Langle \delUb,\, \Big(L_{U_b X} + L_{X U_b}^\top \Big) \delX' \Rangle_S  + \half \Langle \delUd, L_{U_d U_d} \,\delUd' \Rangle \\
    &+ \half \Langle \delUb,\, L_{U_b U_b} \,\delUb' \Rangle_S + \Langle  \delX, F_X^* V_X \Rangle - \Langle V_X, \Delta N \Rangle_S - \Langle V_X, N_{U_b}^\top \delUb \Rangle_S - \Langle \delX, N_X^* V_X \Rangle_S \\
    &+ \Langle V_X,\, F_{U_d}^\top \delUd \Rangle + \Langle \delX ,\,F_X^*  V_{XX} \delX' \Rangle - \Langle  \delX,\,  V_{XX} \Delta N \Rangle_S - \Langle \delX,\,  V_{XX} N_{U_b}^\top \delUb \Rangle_S  \\
    & - \Langle \delX,\, N_X^* V_{XX}\delX' \Rangle_S + \Langle \delX,\, V_{XX} F_{U_d}^\top \delUd \Rangle \bigg].
\end{eqsp}

The form of the \ac{HJB} equation in \cref{eq:HJB_LHS_RHS_green} now properly incorporates boundary information of the value functional. As shown in the subsequent section, the resulting optimal update $\delUb^*$ leverages the first and second derivative of the value functional, which is expected in the context of the established \ac{DDP} method in finite dimensions.

We note that the form of the \ac{HJB} \cref{eq:HJB_LHS_RHS_green} is remarkably different than that of \cite{tzafestas1969differential}. The fundamental differences arise due to a) improper application of Green's theorem, as discussed in \cite{sakawa1972matrix}, and b) terms that are a result of a fundamental difference of reasoning followed in their derivation. For example, the expansions in \cite{tzafestas1969differential} are quite different than the ones computed here, and may reflect an evolution in the \ac{DDP} approach over decades of research.

%===============================================================================

\section{Optimal Distributed and Boundary Control Solutions}

We find two singleton Newton solutions to the \ac{HJB} \cref{eq:HJB_LHS_RHS_green}; one for the optimal distributed control update $\delUd^*$, and one for the optimal boundary control update $\delUb^*$. 

\begin{theorem}
Under the stated assumptions, the optimal distributed update $\delUd^*$ and the optimal boundary update $\delUb^*$ are given in Hilbert spaces by
\begin{align}
    \delUd^* &= - L_{U_d U_d}^{-1}\Big(F_{U_d}^\top V_X + L_{U_d} \Big)  - \half L_{U_d U_d}^{-1}\Big( L_{U_d X} + L_{X U_d}^\top + 2 F_{U_d}^\top V_{XX}\Big)\delX \label{eq:delud} \\
    \delUb^* &= - L_{U_b U_b}^{-1}\Big( L_{U_b} - N_{U_b}^\top V_X \Big) - \half L_{U_b U_b}^{-1} \Big(L_{U_b X} + L_{X U_b}^\top - 2 N_{U_b}^\top V_{XX} \Big) \delX \label{eq:delub}
\end{align}
where we have defined the inverse operators $L_{U_d U_d}^{-1} \in \calL\big(L_2^k(D), L_2^k(D)\big)$ and $L_{U_b U_b}^{-1} \in \calL\big(L_2^l(S), L_2^l(S)\big)$ by their respective inverse kernels, given by
\begin{align}
    L_{U_d U_d}^{-1}(t) W(t) &= \int_D \bar{L}_{U_d U_d}(t, x,y ) W(t, y) \rd y \\
    L_{U_b U_b}^{-1} W(t) &= \int_S \bar{L}_{U_b U_b}(t,\xi,\eta) W(t, \eta) \dSeta
\end{align}
with $\bar{L}_{U_d U_d}(t,x,y)$ and $\bar{L}_{U_b U_b}(t,\xi,\eta)$ denoting the kernel function of the operator $L_{U_d U_d}^{-1}(t)$ and $L_{U_d U_d}^{-1}(t)$ (resp.), and satisfying the property of inverses for kernels
\begin{align}
    \int_D L_{U_d U_d}(t,x,y) \bar{L}_{U_d U_d}\big(t,y,x') \rd y &= I \delta(x-x') \label{eq:d_inverse_property}\\
    \int_S L_{U_b U_b}(t,\xi,\eta) \bar{L}_{U_d U_d}(t,\eta,\xi') \dSeta &= I \delta(\xi-\xi') \label{eq:b_inverse_property}
\end{align}
\end{theorem}

\vspace{1em}
\begin{proof}
    The result can be found by applying a Newton step (e.g. taking the respective partial derivative and setting equal to zero) of the \ac{HJB} \cref{eq:HJB_LHS_RHS_green}. 
    % Details of these steps can be found in the Supplemental Material \textbf{specify section}.
\end{proof}

The equivalent expressions in fields notation expose the spatial integration that takes place in these calculations, and are provided for completeness
\begin{align}
    \delUd &= - \int_D \bar{L}_{U_d U_d}(t,x,y) \Big( F_{U_d}^\top(t,y,y') V_X(t,y') + L_{U_d}(t,y) \Big) \rd y \nonumber \\
    &\quad - \half \int_D \int_D \bar{L}_{U_d U_d}(t,x,y)\Big( L_{U_d X}(t,y,y') + L_{X U_d}^\top(t,y,y') + 2 F_{U_d}(t,y,y'') V_{XX}(t,y'',y')\Big)\delX(t,y') \rd y \; \rd y' \label{eq:delud_fields}\\
    \delUb &= -\int_S \bar{L}_{U_b U_b}(t,x,\eta) \Big(L_{U_b}(t,\eta) - N_{U_b}^\top(t,\eta,\eta') V_X(t, \eta')  \Big) \dSeta \nonumber \\
    &\quad - \half \int_S \int_D \bar{L}_{U_b U_b}(t,x,\eta) \Big( L_{U_b X}(t,\eta, y) + L_{X U_b}(t,\eta,y) - 2 N_{U_b}^\top(t,\eta,\eta') V_{XX}(t,\eta',y)  \Big) \delX(t,y) \rd y \; \dSeta \label{eq:delub_fields}
\end{align}

%===============================================================================

\section{The Backward Value Functional Equations}

The value functional is a backward equation according to the \ac{HJB} \cref{eq:HJB}, and is separated by order into zeroth, first, and second order derivative of the value functional. 
% These backward equations are obtained by plugging in the optimal control updates in \cref{eq:delud,eq:delub} (or \cref{eq:delud_fields,eq:delub_fields} in fields representation) into the \ac{HJB} equation in \cref{eq:HJB_LHS_RHS_green}.
We present the fields representations of these backward equations without cross terms for simplicity. The more general forms of the backward equations with cross terms have been derived, but are lengthy and are omitted due to length considerations.

\begin{theorem}
Under the above stated assumptions, and with optimal control in fields representation given by \cref{eq:delud_fields,eq:delub_fields}, the zeroth-order  backward value functional equation is given by
\small\begin{equation}\label{eq:zeroth_riccati}
\begin{split}
    %%%%%%%%%%%%%%%%ZERO ORDER
    -&\ddt V(t,X(t,x)) \\
    &= L - \int_S V_X(t,\xi)^\top \Delta N(t,\xi) \dSxi \\
    &\quad- \half \int_D \int_D \Big( L_{U_d}^\top(t,y) +  V_X^\top(t,x) F_{U_d}(t,x,y) \Big)\bar{L}_{U_d U_d}(t,y,y') \Big( L_{U_d}(t,y')  + F_{U_d}^\top(t,y',y'') V_X(t,y'') \Big) \rd y\, \rd y' \\
    &\quad - \half \int_S \int_S \Big( L_{U_b}^\top(t,\xi') -  V_X^\top(t,\xi) N_{U_b}(t,\xi,\xi') \Big) \bar{L}_{U_b U_b}(t,\xi',\eta) \Big( L_{U_b}^\top(t,\eta) - N_{U_b}^\top(t,\eta,\eta') V_X(t,\eta') \Big) \rd S_{\xi'}\, \dSeta 
\end{split}
\end{equation}\normalsize
with terminal condition
\begin{equation} \label{eq:zeroth_riccati_IC}
    V\big(t_f,X(t_f,x) \big) = \phi\big(t_f,X(t_f,x) \big),
\end{equation}
the first-order backward value functional equation is given by
\begin{equation} \label{eq:first_riccati}
\begin{split}
    -&\ddt V_X\big(t,X(t,x)\big) \\
    &= L_X(t,x) + F_X^*(t,x,y)V_X(t,y) \\
    &\quad - \int_D \int_D V_{XX}(t,x,x')F_{U_d}(t,x',y) \bar{L}_{U_d U_d}(t,y,y') \Big( L_{U_d}(t,y') + F_{U_d}^\top(t,y',y'') V_X(t,y'') \Big) \rd y \, \rd y' \\
    % boundary integrals
    &\quad + \int_S \int_S V_{XX}(t,\xi,\xi')N_{U_b}(t,\xi',\eta) \bar{L}_{U_b U_b}(t,\eta,\eta') \Big( L_{U_b}(t,\eta') - N_{U_b}^\top(t,\eta',\eta'') V_X(t,\eta'') \Big) \dSeta \, \dSeta' \\
\end{split}
\end{equation}
with boundary and terminal conditions
\begin{align}
    0 &= N_X^*(t,\xi,\eta ) V_X(t,\eta )- \int_S V_{XX}(t,\xi,\eta) \Delta N(t,\eta) \dSeta \label{eq:first_riccati_BC} \\
    &V_X\big(t_f,X(t_f,x)\big) = \phi_X\big(t_f,X(t_f,x)\big), \label{eq:first_riccati_IC}
\end{align}
and the second-order backward value functional equation is given by
\begin{equation}\label{eq:second_riccati}
\begin{split}
    -\ddt V_{XX}(t,x,y) &= L_{XX}(t,x,y) + F_X^*(t,x,y')V_{XX}(t,y',y) + \big[F_X^*(t,x,y') V_{XX}(t,y',y)\big]^\top \\
    &\quad - \int_D \int_D V_{XX}(t,x,x')F_{U_d}(t,x',y')\bar{L}_{U_d U_d}(t,y',y'') F_{U_d}^\top(t,y'',z) V_{XX}(t,z,y)\rd y' \rd y'' \\
    % Boundary terms
    &\quad - \int_S \int_S V_{XX}(t,x,\xi)N_{U_b}(t,\xi,\xi')\bar{L}_{U_b U_b}(t,\xi',\eta) N_{U_b}^\top(t,\eta,\eta') V_{XX}(t,\eta',y)\dSxi' \,\dSeta, 
\end{split}
\end{equation}
with boundary and terminal conditions
\begin{align}
    0 &= N_X^*(t,\xi,\eta) V_{XX}(t, \eta, y) \label{eq:second_riccati_BC} \\
    & V_{XX}\big(t_f,X(t_f,x)\big) = \phi_{XX}\big(t_f, X(t_f, x) \big) \label{eq:second_riccati_IC}
\end{align}
\end{theorem}
\vspace{1em}
\begin{proof}
These equations are obtained by plugging in \cref{eq:delud,eq:delub} into the \cref{eq:HJB_LHS_RHS_green} and grouping terms by order of $\delX$. 
% The intermediate steps can be found in the Supplemental Material \textbf{specify section}.
\end{proof}

The iterative forward-backward system is completed by the approximate variation of the field and boundary dynamics, which are found by rearranging \cref{eq:F_expand,eq:N_expand} as
\begin{align}
    \frac{\rd \delta X(t,x)}{\rd t} &= F(\Xbar + \delX,\Udbar + \delUd,t,x) - F(\Xbar, \Udbar,t,x) \nonumber \\
    &= F_X^\top(\Xbar, \Udbar,t,x) \delX  + F_{U_d}^\top(\Xbar, \Udbar,t,x) \delUd, \quad x \in D \label{eq:delta_dynamics}\\
    0 &= N(\Xbar + \delX, \Ubbar + \delUb, t,\xi) - N(\Xbar, \Ubbar, t, \xi) \nonumber \\
    &= N_X^\top(\Xbar, \Ubbar, t,\xi) \delX  + N_{U_b}^\top(\Xbar, \Ubbar, t,\xi) \delUb, \quad \xi \in S \label{eq:delta_boundary}.
\end{align}
Finally, the control updates of iteration $k+1$ are given by
\begin{align}
    U_d^{k+1} = U_d^k + \gamma_d \delUd^k \label{eq:Ub_update}\\
    U_b^{k+1} = U_b^k + \gamma_b \delUb^k. \label{eq:Ud_update}
\end{align}

%===============================================================================

\section{Recovering Standard Results}

The optimal distributed and boundary control and resulting backward value functional equations represent a generalization of a) \ac{DDP} in finite dimensions and b) the \ac{LQR} for \acp{PDE}. These results are standard results in the control literature, and as such it is important to the validity of our approach to clearly demonstrate that these standard results can be recovered from the equations detailed in the previous sections.

\subsection{Differential Dynamic Programming in Finite Dimensions}

We begin by roughly outlining an analogous derivation of \ac{DDP} in finite dimensions. There are many different formulations of \ac{DDP} in finite dimensions. Our approach specifically follows a body of literature that expands the pertinent functionals around a nominal trajectory. Despite having an extra term for a terminal constraint, we refer to \cite{sun2014continuous} as they present a clean derivation that represents a finite dimensional analogue to the derivation in this document. We ignore the terms having to do with the terminal constraint and the terms that come from second order expansions of the dynamics for ease of comparison. Therein they consider a finite dimensional system of the general form
\begin{equation}
    \ddt x = F(x,u,t), \quad x(t_0) = x_0
\end{equation}
The optimization problem is formulated as 
\begin{eqsp}
    V(x_0, t_0) &= \min_u J(x,u) \\
    &= \min_u \bigg[ \phi\big(x(t_f),t_f\big) + \int_{t_0}^{t_f} L(x,u,t) \rd t \bigg]
\end{eqsp}

After applying standard Taylor expansions of the value functional, its first and second derivative, the dynamics, and the cost functional, plugging them into the \ac{HJB} equation and performing Newton minimization, they obtain the optimal control update as
\begin{equation}
    \delta u = - L_{uu}^{-1}\big( L_u + F_u V_x \big) - \half L_{uu}^{-1}\big( L_{u x} + L_{x u}^\top + 2 F_u V_{XX} \big) \delta x
\end{equation}

This is equivalent in form to the optimal distributed and boundary control update in Hilbert spaces given in \cref{eq:delud,eq:delub}. The resulting backward equations of the value functional in \cite{sun2014continuous} are given by
\begin{align}
    - \ddt V &= L - \half k^\top L_{uu} k \label{eq:finite_zero_riccati}\\
    -\ddt V_x &= L_x + F_x V_x - K^\top L_{uu} k \label{eq:finite_first_riccati}\\
    - \ddt V_{XX} &= L_{xx} - K^\top L_{uu} K + V_{xx} F_x^\top + F_x V_{xx} \label{eq:finite_second_riccati}
\end{align}
where $k \in \Rb^k$, and $K\in \Rb^{k \times n}$ are given by
\begin{align}
    k &= - L_{uu}^{-1}\big( L_u + F_u V_x \big) \\
    K &= - \half L_{uu}^{-1}\big( L_{u x} + L_{x u}^\top + 2 F_u V_{XX} \big)
\end{align}

In order to make the same comparison for the backward equations of the zeroth, first, and second-order value functional in fields, we first define the kernel functions $k_d : T \times \Rb^n \times \Rb^n \rightarrow \Rb^k$, $k_b : T \times \Rb^n \times \Rb^n \rightarrow \Rb^l$, $K_d: T \times \Rb^n \times \Rb^n \times \Rb^n \ \rightarrow T \times \Rb^k \times \Rb^n$, and  $K_d: T \times \Rb^n \times \Rb^n \times \Rb^n \ \rightarrow T \times \Rb^l \times \Rb^n$, which are defined analogously to $k$ and $K$ in \cite{sun2014continuous}, and given by
\begin{align}
    k_d(t,x,y) &= -\bar{L}_{U_d U_d}(t,x,y) \Big(L_{U_d}(t,y) + F_{U_d}^\top(t,y,y') V_X(t,y') \Big) \\
    k_b(t,x,\eta) &= - \bar{L}_{U_b U_b}(t,x,\eta) \Big(L_{U_b}(t,\eta) - N_{U_b}^\top(t,\eta,\eta') V_X(t, \eta')  \Big) \\
    K_d(t,x,y,y')  &= - \half \bar{L}_{U_d U_d}(t,x,y)\Big( L_{U_d X}(t,y,y') + L_{X U_d}^\top(t,y,y')  + 2 F_{U_d}(t,y,y'') V_{XX}(t,y',y'')\Big) \\
    K_b(t,\xi,\eta,y) &= - \half \bar{L}_{U_b U_b}(t,x,\eta) \Big( L_{U_b X}(t,\eta, y) + L_{X U_b}(t,\eta,y) - 2 N_{U_b}^\top(t,\eta,\eta') V_{XX}(t,\eta',y)  \Big)
\end{align}

Thus \cref{eq:zeroth_riccati,eq:first_riccati,eq:second_riccati} in fields representation take the form
\begin{equation}
\begin{split}
    -\ddt V\big(t,X(t,x)\big) &= L - \int_S V_X(t,\xi)^\top \Delta N(t,\xi) \dSxi  \\
    &\quad - \half \int_D \int_D \int_D k_d^\top(t,x,y) L_{U_d U_d}(t,y,y'') k_d(t,y'',y') \rd x\, \rd y \,\rd y' \\
    &\quad - \half \int_S \int_S \int_S k_b^\top(t,\xi,\xi') L_{U_b U_b}(t,\xi',\eta') k_b(t,\eta',\eta) \dSxi\, \dSxi' \, \dSeta  \\
\end{split}
\end{equation}
\begin{equation}
\begin{split}
    -\ddt V_X\big(t,X(t,x)\big) &= L_X + F_X^*V_X - \int_D \int_D \int_D K_d^\top(t,x,y,y') L_{U_d U_d}(t,y',z) k_d(t,z,y'') \rd y \,\rd y'\, \rd y''  \\
    &\quad - \int_S \int_S \int_S K_b^\top(t,x,\xi,\eta,) L_{U_b U_b}(t,\eta,\varphi) k_b(t,\varphi,\eta') \dSxi \; \dSeta \; \dSeta'
\end{split}
\end{equation}
\begin{equation}
\begin{split}
    - \ddt V_{XX}(t,x,y) &= L_{XX}(t,x,y) + F_X^*(t,x,y')V_{XX}(t,y',y) + \big[F_X^*(t,x,y') V_{XX}(t,y',y)\big]^\top  \\
    &\quad - \int_D \int_D \int_D K_d(t,x,y,y')^\top L_{U_d U_d}(t,y',z) K_d(t,z,y'',y''') \rd y' \;\rd y''\; \rd y''' \\
    &\quad - \int_S \int_S \int_S K_b(t,x,\xi, \eta)^\top L_{U_b U_b}(t,\eta,\varphi) K_b(t,\varphi,\eta', y) \dSxi \dSeta \dSeta'
\end{split}
\end{equation}
where in each equation, one of the integrals cancels due to the inverse kernel property in \cref{eq:d_inverse_property,eq:b_inverse_property}. Thus one can recover equations \cref{eq:finite_zero_riccati,eq:finite_first_riccati,eq:finite_second_riccati} by considering an ODE system that a) does not have a spatial state vector so the Volterra-Taylor expansion becomes a Taylor expansion and the volume integrals are equal to their integrand, b) does not have a spatial boundary so surface integrals over the boundary are zero, and c) has real-valued finite dimensional Jacobians defined on an orthonormal basis (with an orthonormal dual basis) so that the adjoint is equal to the transpose.

Thus, the \ac{DDP} equations for \acp{PDE} are a generalization of the \ac{DDP} equations for finite ODE systems. In the following section we demonstrate a similar generalization of the \ac{LQR} solution for \acp{PDE}.

\subsection{The Linear Quadratic Regulator of Fields}

The linear quadratic regulator equations are obtained in \cite{sakawa1972matrix}. Therein, they consider a linear \ac{PDE} of the form
\begin{align}
    \partial_t X(t,x) &= A_x(t) X(t,x) +  B_d(t,x) U_d(t,x), \quad x \in D \\
    X(t_0,x) &= X_0
\end{align}
where $A_x$ is a linear differential operator that has standard form
\begin{equation}
    A_x(t) = \sum_{i,j=1}^n A_{ij}(t,x) \frac{\partial^2}{\partial x_i \partial x_j} + \sum_{i=1}^n B_i(t,x) \frac{\partial}{\partial x_i} + C(t,x).
\end{equation}
The boundary condition is given by
\begin{equation}
    B_b(t,\xi) U_b(t,\xi) = F(t,\xi) X(t,\xi) + \sum_{j=1}^n A_j(t,\xi) \frac{\partial X}{\partial x_j}, \quad \xi \in S
\end{equation}
where the operator $A_j$ is given by
\begin{equation}
    A_j(t,\xi) = \sum_{i=1}A_{ij}(t,\xi) \cos(n_\xi, x_i),
\end{equation}
where $(n_\xi, x_i)$ is the angle between the outward normal $n_\xi$ at a boundary point $\xi \in S$ and the $x_i$-axis. The dynamics are equivalently described in Hilbert spaces as
\begin{align}
    \ddt X(t) &= A(t) X(t) + B_d(t) U_d(t), \quad X \in L_2^n(D) \\
    B_b(t)U_b(t) &= F(t)X(t) + A(t) \cdot \nabla_x X(t), \quad X \in L_2^n(S)
\end{align}
which is a familiar control affine linear system form in Hilbert spaces. The optimization problem is formulated as
\begin{equation}
    J(t_0, X_0, U_d, U_b) = \phi\big(X(t_f)\big) + \int_{t_0}^{t_f} L\big(t,X(t), U_d(t), U_b(t) \big) \rd t
\end{equation}
where the running cost $L$ has the form
\begin{eqsp}
    L\big(t,X(t),U_d(t), U_b(t) \big) &= \half \int_D \int_D X(t,x)^\top Q(t,x,y) X(t,y) \rd x \rd y + \half \int_D \int_D U_d(t,x)^\top R_d(t,x,y) U_d(t,y) \rd x \rd y \\
    &\quad + \half \int_S \int_S U_b(t,\xi)^\top R_b(t,\xi,\eta) U_b(t,\eta) \dSxi \dSeta
\end{eqsp}
where the kernels $Q \geq 0$, $R_d>0$ and $R_b>0$ are all assumed to be symmetric about all spatial axes.

% Just as above, the \ac{HJB} is given by

The resulting optimal distributed and boundary control equations are obtained after applying Green's theorem to the \ac{HJB} equation and performing Newton minimization. They are given by
\begin{align}
    U_d^*(t) &= - R_d^{-1}(t) B_d(t)^\top P(t) X(t) \label{eq:LQR_Ud}\\
    U_b^*(t) &= - R_b^{-1}(t) B_b(t)^\top P(t) X(t)\label{eq:LQR_Ub}
\end{align}
which are equivalently written in fields representation as
\begin{align}
    U_d^*(t,x) &= - \int_D \int_D \bar{R}_d(t,x,y) B_d^\top(t,y) P(t,y,y') X(t,y') \rd y \rd y' \label{eq:LQR_Ud_fields} \\
    U_b^*(t,x) &= -\int_S \int_D \bar{R}_b(t,\xi,\eta) B_b^\top(t,\eta) P(t,\eta,y) X(t,y) \rd y \dSeta \label{eq:LQR_Ub_fields}
\end{align}

\begin{table}[t!]
\centering
\captionof{table}{\label{table:LQR_DDP_equivalences}Corresponding Hilbert space operators between \ac{LQR} of fields and \ac{DDP} of fields.}
\begin{tabular}{@{} c c @{}}
    \toprule
    \textbf{DDP Operators} & \textbf{LQR Operators} \\
    \midrule
    $L_{U_d U_d}(t)$ & $R_d(t)$\\
    $L_{U_b U_b}(t)$ & $R_b(t)$\\
    $L_{XX}(t)$ & $Q(t)$ \\
    $F_{U_d}(t)$ & $B_d(t)$\\
    $-N_{U_b}(t)$ & $B_b(t)$\\
    $V_{X}(t)$ & $P(t)X(t)$ \\
    $V_{XX}(t)$ & $P(t)$ \\
    \bottomrule
\end{tabular}
\end{table}

%%%%%%%%%%%%%%%%%OLD VERSION OF THE TABLE%%%%%%%%%%%%
% \begin{table}[t!]
% \centering
% \captionof{table}{\label{table:LQR_DDP_equivalences}Corresponding Operators between \ac{LQR} of fields and \ac{DDP} of fields. All operators are written in Hilbert spaces for simplicity. \oswin{Maybe consider using booktabs for the table?}}
% \begin{tabular}{c|c}
%     \textbf{DDP Operators} & \textbf{LQR Operators} \\
%     \hline
%     $L_{U_d U_d}(t)$ & $R_d(t)$\\
%     $L_{U_b U_b}(t)$ & $R_b(t)$\\
%     $L_{XX}(t)$ & $Q(t)$ \\
%     $F_{U_d}(t)$ & $B_d(t)$\\
%     $-N_{U_b}(t)$ & $B_b(t)$\\
%     $V_{X}(t)$ & $P(t)X(t)$ \\
%     $V_{XX}(t)$ & $P(t)$ \\
% \end{tabular}
% \end{table}
%%%%%%%%%%%%%%%%%%%%%%%%%%%%%%%%%%%%%%%%%%%%%%%%%%%%%

In order to make the generalization clear, we rewrite \cref{eq:LQR_Ud,eq:LQR_Ub} in our notation using the correspondences listed in \cref{table:LQR_DDP_equivalences}
\begin{align}
    U_d^*(t) &= - L_{U_d U_d}^{-1}(t) F_{U_d}^\top(t) V_X(t,X) \label{eq:LQR_Ud_DDP_notaion}\\
    U_b^*(t) &= L_{U_b U_b}^{-1}(t) N_{U_b}^\top(t) V_X(t,X) \label{eq:LQR_Ub_DDP_notaion}
\end{align}
and repeat \cref{eq:delud,eq:delub} here for clarity
\begin{align*}
    \delUd^* &= \textcolor{blue}{- L_{U_d U_d}^{-1}} \Big(  L_{U_d} + \textcolor{blue}{F_{U_d}^\top V_X} \Big) - \half L_{U_d U_d}^{-1}\Big( L_{U_d X} + L_{X U_d}^\top + 2 F_{U_d}^\top V_{XX}\Big)\delX  \\
    \delUb^* &= \textcolor{blue}{- L_{U_b U_b}^{-1}}\Big( L_{U_b} \textcolor{blue}{- N_{U_b}^\top V_X} \Big) - \half L_{U_b U_b}^{-1} \Big(L_{U_b X} + L_{X U_b}^\top - 2 N_{U_b}^\top V_{XX} \Big) \delX 
\end{align*}
Thus, we can recover \cref{eq:LQR_Ud_DDP_notaion,eq:LQR_Ub_DDP_notaion} by a) assuming the cost functional is a pure quadratic without cross terms so that the terms $L_{U_d}$, $L_{U_b}$, $L_{X U_d}, L_{U_d X}, L_{X U_b}$, and $L_{U_b X}$ are null and b) using only gradient information of the value functional so that $V_{XX}$ terms are ignored.

The resulting second-order backward value functional equation of Riccati type for \ac{LQR} is given in fields representation by
\begin{equation}
    \begin{split}
    \frac{\partial P(t,x,y)}{\partial t} &= - A_X^*(t) P(t,x,y) - \big[ A_X^*(t) P(t,x,y) \big]^\top  - Q(t,x,y) \\
    &\quad + \int_D \int_D P(t,x,x') B_d(t,x') \bar{R}_d(t,x',x'') B_d^\top(t,x'')P(t,x'',y) \rd x' \rd x'' \\
    &\quad + \int_S\int_S P(t,x,\xi) B_b(t,\xi) \bar{R}_b(t,\xi,\eta) B_b^\top(t,\eta) P(t,\eta,y) \dSxi \dSeta
    \end{split}
\end{equation}
which is rewritten in the DDP notation by again applying the correspondences listed in \cref{table:LQR_DDP_equivalences} as
\begin{equation} \label{eq:LQR_riccati_DDP_notation}
\begin{split}
    \frac{\partial V_{XX}(t,x,y)}{\partial t} &= - F_X^*(t,x,y) V_{XX}(t,x,y) - \big[ F_X^*(t,x,y) V_{XX}(t,x,y) \big]^\top - L_{XX}(t,x,y) \\
    &\quad + \int_D \int_D V_{XX}(t,x,x') F_{U_d}(t,x') \bar{L}_{U_d U_d}(t,x',x'') F_{U_d}^\top(t,x'') V_{XX}(t,x'',y) \rd x' \rd x'' \\
    &\quad + \int_S\int_S V_{XX}(t,x,\xi) F_{U_b}(t,\xi) \bar{L}_{U_b U_b}(t,\xi,\eta) F_{U_b}^\top(t,\eta) V_{XX}(t,\eta,y) \dSxi \dSeta.
\end{split}
\end{equation}
\cref{eq:LQR_riccati_DDP_notation} is identical to the second order backward value functional of \ac{DDP} of fields without cross terms in the running cost, given in \cref{eq:second_riccati}. 

Thus we conclude that the equations of \ac{STDDP} are a generalization of \ac{LQR} of fields. This generalization is analogous to the similar generalization of \ac{LQR} of ODE systems to \ac{DDP} of ODE systems. Whereas \ac{LQR} is the analytically optimal controller for linear systems, it cannot be applied directly to a nonlinear system, nor can it be applied directly to a linear system whose running cost functional is not purely quadratic. In contrast, the iterative approximate optimal control method provided by \ac{DDP} of fields was constructed for such systems.

%===============================================================================

\section{Continuous-Time Convergence Analysis}

Global convergence of the discrete finite dimensional \ac{DDP} algorithm defined for discrete ODE systems was first provided by Yakowitz and Rutherford \cite{yakowitz1984computational}. Later the proof that discrete finite dimensional \ac{DDP} converges quadratically in the number of iterations was proved independently by Pantoja \cite{pantoja1983algorithms} and Murray and Yakowitz \cite{murray1984differential}. This quadratic convergence proof relied on convergence of Newton's method, but later an independent proof relying only on the dynamic programming principle was given by Liao and Shoemaker \cite{shoemaker1990proof}. Through decades of application of the \ac{DDP} algorithm, there have been numerous extensions of the proof of global convergence, for example for \ac{DDP} on Lie groups in \cite{boutselis2018differential} and for \ac{DDP} with generalized Polynomail Chaos expansions in \cite{boutselis2019numerical}. However it appears to the best knowledge of the authors that most if not all proofs of global convergence are for \ac{DDP} and its extensions in discrete time, and not in continuous time. 

Typically, one determines provable convergence characteristics by investigating the behavior of the derivative
\begin{equation}
    \frac{\rd J^i\big(t, X(t), U(t)\big)}{\rd i} = \frac{\rd J^i\big(t, X(t), U(t)\big) }{ \rd U^i_{t_0:t_f}} \frac{\rd U^i_{t_0:t_f} }{\rd i}
\end{equation}
where, due to the decoupled nature of the distributed and boundary control updates, we have defined the Hilbert space control vector $U(t) \in  L_2^{k+l}(\bar{D})$ as the direct product Hilbert space analog of the stacked distributed and boundary control vectors in fields representation $U(t,x) = [U_d(t,x),\, U_b(t,x)]^\top$. This notation simplifies our analysis significantly. We have also introduced the trajectory notation, where subscript $t_0\!:\!t_f$ represents the entire trajectory in time of the associated variable, and used the superscript $i$ for the \ac{STDDP} iteration index. Similarly $\delU_{t_0:t_f}^i$ denotes the control update trajectory for control trajectory $U_{t_0:t_f}^i$. This trajectory notation defines a \textit{temporal} Hilbert space over time-indexed \textit{spatial} Hilbert spaces. Let $L_2^n(T)$ denote the Hilbert space of $L_2^n(\bar{D})$-vector functions square integrable over $T$ with inner product
\begin{equation}
    \Langle X_{1,t_0:t_f}, X_{2,t_0:t_f} \Rangle_T = \int_{t_0}^{t_f} \Langle X_{1}(s), X_{2}(s) \Rangle \rd s.
\end{equation}
This allows us to write time integrals over trajectory variables as inner product tensor contractions, and treat continuous trajectories as objects in a similar way to the continuum of the \ac{PDE} variables. We begin by stating the following lemma, assumption, and proposition that will be used in our analysis.

\begin{lemma}
Assume the cost functional has the form of \cref{eq:hilbert_cost} and define the measurable backward recursive functional 
% $\psi:T\times L_2^n(\bar{D}) \times L_2^{k+l}(\bar{D}) \rightarrow \calL\big(L_2^n(\bar{D})\big)$ for some $\varepsilon >0$ as
$\psi \in \calL\big(L_2^n(\bar{D})\big)$ for some $\varepsilon > 0$ as
\begin{align}
    \psi\big(t,X(t),U(t)\big) &= \int_t^{t+\varepsilon} L_X\big(s,X(s),U(s)\big) \rd s + \Phi(t,s)\psi\big(t + \varepsilon, X(t+\varepsilon), U(t+\varepsilon) \big) \label{eq:psi} \\
    \psi\big(t_f, X(t_f), U(t_f) \big) &= \phi_X\big(t, X(t) \big) \label{eq:psi_BC}
\end{align}
where $\Phi(\cdot, \cdot)$ is a contractive linear semigroup generated by the approximate variation dynamics in \cref{eq:delta_dynamics,eq:delta_boundary}, and is assumed to be positive definite almost everywhere. Then the cost functional satisfies
\begin{equation}\label{eq:conv_lemma}
    \frac{\rd J^i\big(t, X(t), U(t)\big) }{ \rd U^i_{t_0:t_f}} = \Langle L_{U,t_0:t_f}, \mathds{1}_{t_0:t_f} \Rangle_T + \Langle F_{U,t_0:t_f}, \psi_{t_0:t_f} \Rangle_T
\end{equation}
where $\mathds{1}_{t_0:t_f}$ is the trajectory of ones.
\newline
\begin{proof}
    The proof is in the Supplementary Material, \cref{sec:proof_convergence_lemma}
\end{proof}
\end{lemma}

\begin{assumption}\label{assumption_compactness}
The search space of control trajectories $\calU \ni U_{t_0:t_f}^i$ is compact.
\end{assumption}

Our analysis is simplified by the $Q$ functional notation defined as follows
\begin{align*}
    \begin{array}{l l}
    Q_{UU} = L_{UU} &  Q_U = L_U + F_U^\top V_X \\
    Q_{UX} = L_{UX} + F_U^\top V_{XX} & Q_X = L_X + F_X^\top V_X \\  Q_{XX} = L_{XX} + F_X^*V_{XX} +  \big[ F_X^* V_{XX}\big]^\top & 
    \end{array}
\end{align*}

\begin{prop}\label{D_prop}
Let the \ac{PDE} $D(t) \in L_2^n(\bar{D})$ have dynamics
\begin{align}
    \ddt D(t) &= -F_X^\top D(t) + Q_{UX}^\top Q_{UU}^{-1} Q_U \label{eq:D_dynamics}\\
    D(T) &= 0 \label{eq:D_boundary}
\end{align}
Then $D(t)$ has weak backwards solutions defined in the Hadamard sense, and given by
\begin{equation}\label{eq:D_sol}
    D(t) = \int_T^t \Phi^\top(t,\tau) Q_{UX}^\top(\tau) Q_{UU}^{-1}Q_U(\tau) \rd \tau
\end{equation}
\begin{proof}
    The existence of weak solutions is given by the assumption that solutions to $\psi$ and $V_X$ exist. The rest of the proof is immediate given that the dynamics are of semilinear form and have a zero terminal condition.
\end{proof}
\end{prop}

\begin{theorem} \label{thm:quadratic}
Consider the continuous-time optimal control problem in \cref{eq:val_fcn} subject to the dynamics in \cref{eq:hilbert_dynamics,eq:hilbert_boundary} 
with cost functional $J$ having no cross terms for simplicity. 
Let $\bar{U}_{t_0:t_f} \in L_2^{k+l}(T)$ be a nominal control trajectory and let $\delU_{t_0:t_f} \in L_2^{k+l}(T)$ be the trajectory of control updates from \cref{eq:delud,eq:delub}. Then the following holds:
\begin{equation}\label{eq:conv_thm}
    \frac{\rd J^i\big(t, X(t), U(t)\big)}{\rd i} = -\gamma \Langle Q_{U,t_0:t_f}, M_{t_0:t_f} Q_{U,t_0:t_f} \Rangle_T + O(\gamma^2)
\end{equation}
where the trajectory operator $M_{t_0:t_f} \in \calL\big(L_2^{l+k}(T), L_2^{l+k}(T) \big)$ has a positive definite kernel $\forall t \in T$.
\newline
\begin{proof}
% Assume for simplicity that $J$ has no cross terms, i.e. $L_{UX}=0$, and
Observe that due to the iterative updates in \cref{eq:Ud_update,eq:Ub_update}, we have
\begin{equation}
    \frac{\rd U_{t_0:t_f}}{\rd i} = \gamma \delU_{t_0:t_f}.
\end{equation}
Also at our disposal is the identity
\begin{equation} \label{eq:identity}
    L_{U,t_0:t_f} + F_{U,t_0:t_f}^\top \psi_{t_0:t_f} = Q_{U,t_0:t_f} - F_{U,t_0:t_f}^\top \big(V_{X,t_0:t_f} - \psi_{t_0:t_f} \big).
\end{equation}
Plugging \cref{eq:identity} into \cref{eq:conv_lemma} yields
\begin{align}
    \frac{\rd J^i\big(t, X(t), U(t)\big)}{\rd i}
    &= -\gamma \Langle Q_{U,t_0:t_f}, Q_{UU,t_0:t_f}^{-1} Q_{U,t_0:t_f} \Rangle_T + \gamma \Langle F_{U,t_0:t_f}^\top \big(V_{X,t_0:t_f} - \psi_{t_0:t_f} \big), Q_{UU,t_0:t_f}^{-1} Q_{U,t_0:t_f} \Rangle_T \nonumber \\
    &\quad -\gamma \Langle Q_{U,t_0:t_f}, Q_{UU,t_0:t_f}^{-1} Q_{UX,t_0:t_f} \delX_{t_0:t_f} \Rangle_T \nonumber \\
    &\quad + \gamma \Langle F_{U,t_0:t_f}^\top \big(V_{X,t_0:t_f} - \psi_{t_0:t_f} \big), Q_{UU,t_0:t_f}^{-1} Q_{UX,t_0:t_f} \delX_{t_0:t_f} \Rangle_T
\end{align}
Note that the total variation $\delX \in L_2^n(\bar{D})$, with dynamics given in semilinear form by \cref{eq:delta_dynamics,eq:delta_boundary}, has a solution given by
\begin{equation}
    \delX(t) = \Phi(t,t_0)\delX_0 + \int_0^t \Phi(t,s)F_U(s,)\delU(s) \rd s
\end{equation}
Thus, since $\delX_0 = 0$, and $\delU(t) = O(\gamma)$, it follows that $\delX = O(\gamma)$, so we have
\begin{align}
    \frac{\rd J^i\big(t, X(t), U(t)\big)}{\rd i} &= -\gamma \Langle Q_{U,t_0:t_f}, Q_{UU,t_0:t_f}^{-1} Q_{U,t_0:t_f} \Rangle_T  \nonumber \\
    &\quad + \gamma \Langle F_{U,t_0:t_f}^\top \big(V_{X,t_0:t_f} - \psi_{t_0:t_f} \big), Q_{UU,t_0:t_f}^{-1} Q_{U,t_0:t_f} \Rangle_T + O(\gamma^2) \label{eq:dJ_di}
\end{align}
Next, notice that $D(t) = V_{X}(t) - \psi(t)$ has dynamics of the form of \cref{eq:D_dynamics}, with an equivalent terminal condition. Thus, plugging in \cref{eq:D_sol} into \cref{eq:dJ_di} and reducing yields
\begin{align}
    \frac{\rd J^i\big(t, X(t), U(t)\big)}{\rd i} &= -\gamma \Langle Q_{U,t_0:t_f}, Q_{UU,t_0:t_f}^{-1} Q_{U,t_0:t_f} \Rangle_T - \gamma \Langle  Q_{U,t_0:t_f}, A_{1,t_0:t_f} Q_{U,t_0:t_f} \Rangle_T + O(\gamma^2)  \label{eq:dJ_di_reduced}
\end{align}
where $A_{1,t_0:t_f} \in \calL\big(L_2^{k+l}(T),L_2^{k+l}(T)\big)$ has a positive definite kernel $\forall$ $t \in T$ due to the positive definiteness of $\Phi(\cdot,\cdot)$ by definition, and the positive definiteness of $V_{XX}$ by assumption. Thus, due to the positive definiteness of the kernels of $L_{U_d U_d}(t)$ and $L_{U_b U_b}(t)$ by assumption, one can form $M \in \calL\big(L_2^{k+l}(\bar{D}),L_2^{k+l}(\bar{D})\big)$ with positive definite kernel $\forall$ $t \in T$ such that
\begin{equation}
    \frac{\rd J^i\big(t, X(t), U(t)\big)}{\rd i} = -\gamma \Langle Q_{U,t_0:t_f}, M_{t_0:t_f} Q_{U,t_0:t_f} \Rangle_T + O(\gamma^2)
\end{equation}
which concludes the proof
\end{proof}
\end{theorem}

\begin{corollary}
Suppose \cref{assumption_compactness} holds. Then the iterative  \cref{eq:hilbert_dynamics,eq:hilbert_boundary,eq:delud,eq:delub,eq:zeroth_riccati,eq:zeroth_riccati_IC,eq:first_riccati,eq:first_riccati_IC,eq:first_riccati_BC,eq:second_riccati,eq:second_riccati_BC,eq:second_riccati_IC,eq:Ud_update,eq:Ub_update} will converge to a stationary solution of \cref{eq:val_fcn}.
\newline
\begin{proof}
    The proof is in the Supplementary Material,  \cref{sec:proof_convergence_corollary}
\end{proof}
\end{corollary}

% My plan for the proof:
% \begin{itemize}
%     \item In the paper: theorem that states that $\frac{\rd J}{\rd i} = -\gamma Z^\top(t) M(t) Z(t)$, $\forall$ $t \in T$
%     \item in the paper: theorem that states that the iterative approach provided by \cref{Algorithm1} without cost functional cross terms converges a.e. to the state and control triple $(X^*, U_d^*, U_b^*)$ for any set of initial conditions $X_0$ and control initialization $\Udbar$, $\Ubbar$
%     \item in supplement: remark that since the optimal control updates $\delUd$, $\delUb$ are decoupled, we can define a new Hilbert space $\calD$ and have a single dynamical equation (i.e. the deterministic form of \cite{fabbri}, with $\hat{X} \in \calD$, $\hat{U} \in \calD$ and analogous linearization terms $\hat{F}_X$, $\hat{F}_U$.
%     \item in supplement: corollary that the \ac{PDE} solution can be approximated by the semi-group of the linearization.
%     \item With this change of notation, proof of the lemma/theorem follow from the overleaf document
% \end{itemize}

%===============================================================================

\section{STDDP Algorithm}

The resulting \ac{STDDP} algorithm can be applied for control of any nonlinear forward spatio-temporal \ac{PDE} system satisfying the stated assumptions. It is an iterative forward-backward approach, wherein each iteration forward propagates the dynamics, backward propagates the value functional and its derivatives, and updates the control based on approximate variation dynamics. The resulting procedure is described in greater detail in \cref{Algorithm1}. 

Note that \cref{Algorithm1} has a forward process, a backward process, and another forward process. While this is a simpler algorithmic exposition, the runtime performance can be improved by simply combining the two forward time loops. While the numerical experiments in this manuscript were performed with a fixed learning rate for demonstration purposes, it can be numerically advantageous to apply line search methods to adapt the learn rate. Some such methods are described in \cite{tassa2012synthesis} and \cite{tassa2014control}, and typically evaluate the best learning rate based on the best improvement in the cost functional. However, since the value functional typically encodes problem information beyond the cost metric, one may also evaluate learning rate based on improvements in the value functional.

The inputs of the \ac{STDDP} algorithm can change depending on the specific problem but in most cases contain time interval ($T$), number of iterations ($K$), initial state ($X_0$), time discretization ($\Delta t$), distributed control learn rate ($\gamma_d$), and boundary control learn rate ($\gamma_{b}$). One may also include a number of rollouts ($R$) for a parallelized line search. Instead of a fixed number of iterations, one may also check for convergence using relative or absolute convergence criteria in either the cost functional or the value functional \cite{tassa2014control}.

\begin{algorithm}[t]
 \caption{STDDP}
 \begin{algorithmic}[1]
 \State \textbf{Function:} \textit{$(U_d^*, U_b^*) =$ \textbf{STDDP}($T$,$K$,$X_0$, $\Udbar$, $\Ubbar$, $\Delta t$, $\gamma_d$, $\gamma_b$)}
 \For {$k=1 \; \text{to} \; K$}
 \State Forward propagate PDE dynamics in \cref{eq:hilbert_dynamics}
 \State Evaluate running cost $L$ and its partial derivatives
 \State Evaluate terminal cost $\phi$ and its partial derivatives
 \State Backward propagate value functional via \cref{eq:zeroth_riccati,eq:first_riccati,eq:second_riccati}
 \State Forward propagate approximate variation dynamics via \cref{eq:delta_dynamics,eq:delta_boundary}
 \State Compute updates $\delUd^k$ and $\delUb^k$ via \cref{eq:delud,eq:delub}
 \State Update control $U_d^{k+1}$, $U_b^{k+1}$ via \cref{eq:Ud_update,eq:Ub_update}
 \EndFor
 \end{algorithmic}
 \label{Algorithm1}
\end{algorithm}

\subsection{Forward \& Backward PDE Discretization Methods} \label{sec:discretization}

In order to implement the forward spatio-temporal system dynamics in \cref{eq:hilbert_dynamics,eq:hilbert_boundary} and the backward value functional system in \cref{eq:zeroth_riccati,eq:zeroth_riccati_IC,eq:first_riccati,eq:first_riccati_BC,eq:first_riccati_IC,eq:second_riccati,eq:second_riccati_BC,eq:second_riccati_IC} on a digital computer, these forward and backward \acp{PDE} must be spatially and temporally discretized. 

Nonlinear \acp{PDE} in the Eulerian formalism are often spatially discretized using either finite difference methods, Galerkin methods, or finite element methods. In this work we apply a spatial central finite difference discretization, which yields a fixed 1D grid of length $a$, with $J$ elements. We note that through the above derivation, any discretization can be used in place of the central difference.

While there are numerous works describing temporal discretization methods for a multitude of forward \acp{PDE}, there are relatively few that describe temporal discretization schemes for backward \acp{PDE} of Riccati type. In finite dimensions, these are typically referred to as \acp{RDE}, and their discretization presents several difficulties which stem from a matrix-valued variable that cannot be analytically isolated without using a Kronecker scheme. Furthermore, \acp{RDE} are known to be quite stiff in many contexts due to a fast transient response \cite{butusov2016semi}. 

The most straightforward method is the explicit time Euler discretization method, which has a fast implementation, yet is sensitive to discretization step size for stiff dynamics. This sensitivity can be reduced by applying Runge-Kutta time-integration techniques, however one must either super-sample the dynamics or apply an equivalent Runge-Kutta integration for the dynamics and value functionals.

Semi-implicit time discretization and implicit time discretization are well known to handle stiff dynamics, yet require isolation of the value functional. This in turn yields an update with a very large Kronecker sum matrix inversion. 
% To see this more clearly, consider the finite dimensional RDE given by
% \begin{equation}
%     -\dot{V}_{xx}(t) = L_{xx} + F_x^\top V_{xx}(t) + V_{xx}(t) F_x - V_{xx}(t) F_u L_{uu}^{-1} F_u^\top V_{xx}(t)
% \end{equation}
% This can be equivalently written in vector form by stacking columns of matrices with the $\rvec$ operator
% \begin{align}
%     - \rvec\big(\dot{V}_{xx}(t)\big) &= \rvec(L_{xx}) + F_x^\top \oplus F_x^\top \rvec\big(V_{xx}(t)\big) \nonumber \\
%     &\quad - \rvec\big( V_{xx}(t) F_u L_{uu}^{-1} F_u^\top V_{xx}(t) \big)
% \end{align}
% To elucidate, consider the finite dimensional analog of
To elucidate, consider the discretized \ac{1D} Hilbert space representation of \cref{eq:second_riccati}, where $F_X^* = F_X^\top$, given by
\begin{equation}
\begin{split}
    -\ddt V_{XX}(t) &= L_{XX} + \frac{1}{\Delta x} F_{X}^\top V_{XX}(t) + \frac{1}{\Delta x} V_{XX}(t) F_X - \frac{1}{\Delta x^2} V_{XX}(t) F_{U_d} L_{U_d U_d}^{-1} F_{U_d}^\top V_{XX}(t) \\
    &\quad - \frac{1}{\Delta x^2} V_{XX}(t) N_{U_b} L_{U_b U_b}^{-1} N_{U_b}^\top V_{XX}(t).
\end{split}
\end{equation}
% \begin{equation}
% \begin{split}
%     -\ddt V_{xx}(t) &= L_{xx} + F_{x}^\top V_{xx}(t) + V_{xx}(t) F_x - V_{xx}(t) F_{uu} L_{uu}^{-1} F_{u}^\top V_{xx}(t)
% \end{split}
% \end{equation}

Clearly, the desired variable $V_{XX}$ cannot be completely isolated in this form. However, one can equivalently write a vector form by application of the $\rvec$ operator
\begin{equation}
\begin{split}
    - \rvec\left(\ddt V_{XX}(t)\right) &= \rvec(L_{XX}) + \frac{1}{\Delta x} F_X^\top \oplus F_X^\top \rvec\big(V_{XX}(t) \big) - \frac{1}{\Delta x^2} \rvec\big(V_{XX}(t) F_{U_d} L_{U_d U_d}^{-1} F_{U_d}^\top V_{XX}(t)\big) \\
    &\quad - \frac{1}{\Delta x^2} \rvec\big(V_{XX}(t) N_{U_b} L_{U_b U_b}^{-1} N_{U_b}^\top V_{XX}(t)\big).
\end{split}
\end{equation}
% \begin{equation}
% \begin{split}
%     - \rvec\left(\ddt V_{xx}(t)\right) &= \rvec(L_{xx}) + F_x^\top \oplus F_x^\top \rvec\big(V_{xx}(t) \big) \\
%     &\quad - \rvec\big(V_{xx}(t) F_{u} L_{uu}^{-1} F_{u}^\top V_{xx}(t)\big)
% \end{split}
% \end{equation}

Semi-implicit time discretization schemes typically evaluate terms that are linear in $V_{XX}(t)$ at the current time step and non-linear terms in $V_{XX}(t)$ at the next time step \cite{lord_powell_shardlow_2014}, which is the previous time step in the case of backward \acp{PDE}. The resulting semi-implicit update is given by
\begin{equation}
    \begin{split}
    \rvec\Big( &V_{XX}(t_{k-1}) \Big) \\
    &= \bigg[ I - F_X^\top \otimes F_X^\top \Delta t \bigg]^{-1} \bigg[ \rvec\Big(V_{XX}(t_{k})  \Big)  + \rvec\Big( L_{XX}\Big) \Delta t - \rvec\Big( V_{XX}(t_{k})F_{U_d} L_{U_d U_d}^{-1}F_{U_d}^\top V_{XX}(t_{k}) \Big)\Delta t \bigg] .
    \end{split}
\end{equation}
% \begin{equation}
%     \begin{split}
%     &\rvec\Big( V_{xx}(t_{k-1}) \Big) \\
%     &= \bigg[ I - F_x^\top \otimes F_x^\top \Delta t \bigg]^{-1} \bigg[ \rvec\Big(V_{xx}(t_{k})  \Big)  + \rvec\Big( L_{xx}\Big) \Delta t \\
%     &\qquad\qquad\qquad\qquad\;\;\;\, - \rvec\Big( V_{xx}(t_{k})F_{u} L_{uu}^{-1}F_{u}^\top V_{xx}(t_{k}) \Big)\Delta t \bigg] 
%     \end{split}
% \end{equation}
The resulting update equation is less sensitive to time discretization step size $\Delta t$, however it requires the inversion of a large matrix of size $J^2 \times J^2$ for each time step of each iteration, where $J$ is the spatial discretization size of the \ac{1D} \ac{PDE}. A key observation is that the matrix $M := I - F_X^\top \otimes F_X^\top \Delta t$ typically only has as many diagonals as the order of the spatial discretization, and is zeros elsewhere except for the boundary conditions, thus it is a sparse matrix. For example, in the case of a second order spatial central difference discretization of the Burgers equation with Homogeneous Dirichlet boundary conditions, $M$ is tridiagonal. Thus the inverse can be efficiently computed with sparse linear equation solvers such as SuperLU \cite{demmel1999supernodal}.

In \cite{butusov2016semi}, the authors describe so called D-methods, which reduce computational complexity inherent to semi-implicit methods by applying explicit Euler discretization to some subset of the variables, and applies semi-implicit discretization to the rest. This could dramatically reduce complexity; if $J_e \leq J$ is the number of points treated with explicit discretization, then the resulting semi-implicit inverse is of size $(J-J_e)^2 \times (J-J_e)^2$. This may have dramatic benefit for \acp{ODE} systems where one may have slower and faster channels, However it is not clear how to select grid elements for the associated D-method for Riccati \acp{PDE}.

%===============================================================================

\section{Simulated Experiments}

We applied the \ac{STDDP} algorithm to two simulated \ac{PDE} experiments to optimally control the system to a prescribed desired behavior. Each experiment used less than 32 GB RAM, and was run on a desktop computer with an Intel Xeon 12-core CPU with a NVIDIA GeForce GTX 980 GPU. The computations did not utilize GPU parallelization, however many operations, such as cost and partial derivative computations, can be parallelized for greater computational efficiency. 

The simulated experiments involve reaching tasks, where the \ac{PDE} is initialized at a zero initial condition over the spatial region, and must reach certain field values at prescribed regions of the spatial domain. As discussed in the previous section, each \ac{PDE} was spatially discretized by a spatial central difference discretization, and an expicit-time Euler discretization. The first and second derivative of the value functional were spatially discretized on the same spatial central difference grid as the forward dynamical \ac{PDE}, and all three backward equations were temporally discretized with an explicit Euler discretization. Regularization was added to the second derivative of the value functional in order to aid in numerical stability.

%%%%%%%%%HEAT FIGURES%%%%%%%%%%%%%%%%
{\centering
\begin{figure*}[!t]
    % \subfigure[Controlled Contour]{\includegraphics[width=0.52\textwidth]{Figures/1DHeat_contour_tfpt12_scost281456_scost_w300_cropped.png}}
    % \subfigure[Final Time Snapshot]{\hspace{.1cm}\includegraphics[width=0.470\textwidth]{Figures/1DHeat_final_time_tfpt12_scost281456_scost_w300_1_cropped.png}}
    
    \begin{multicols}{2}
    \begin{subfigure}[h!]{1.0\textwidth}
    \hspace{0.0cm}\includegraphics[width=0.52\textwidth]{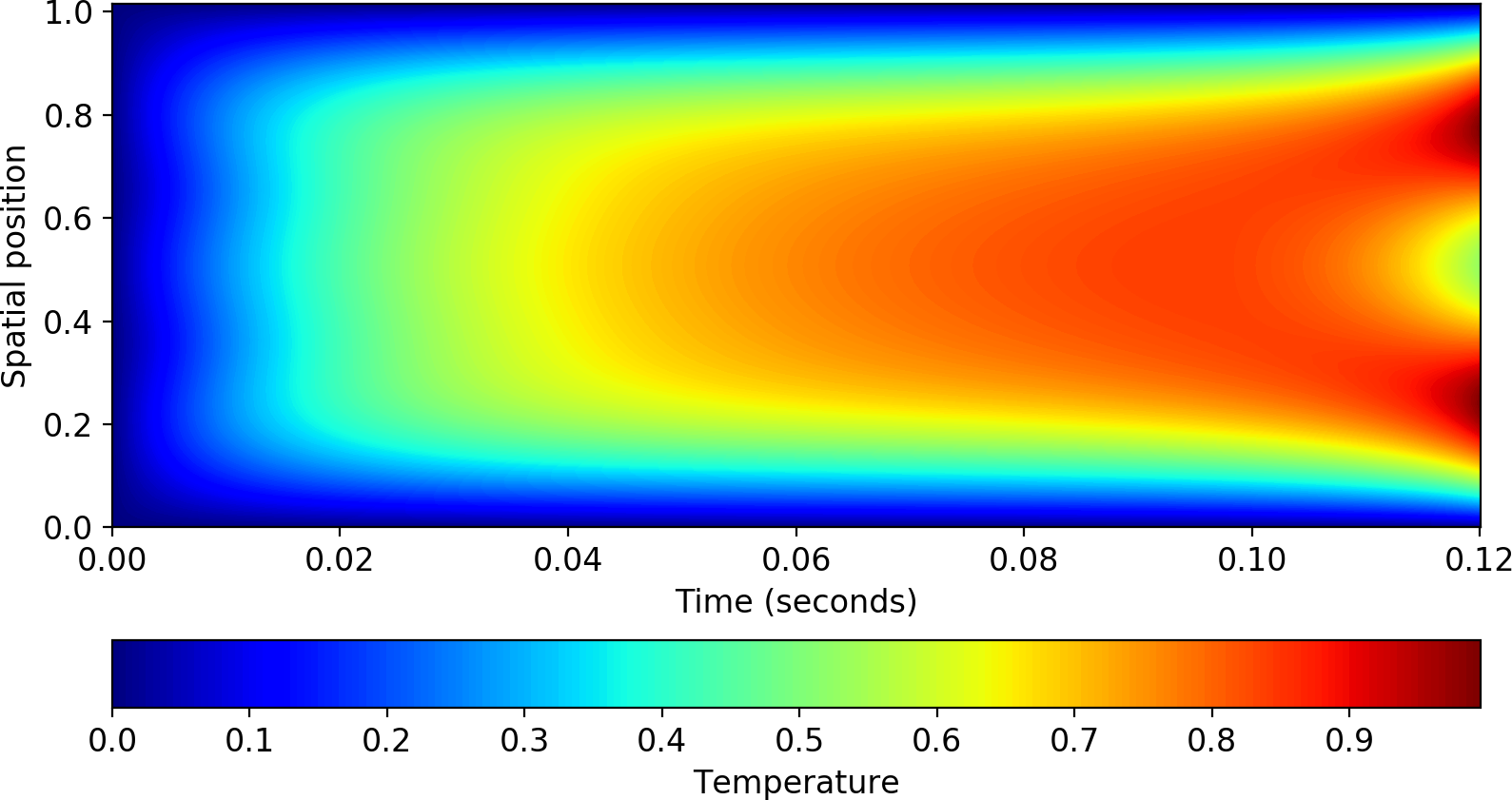}
    \end{subfigure}
    
    \begin{subfigure}[h!]{1.0\textwidth}
     \hspace{.4cm}\includegraphics[width=0.46\textwidth]{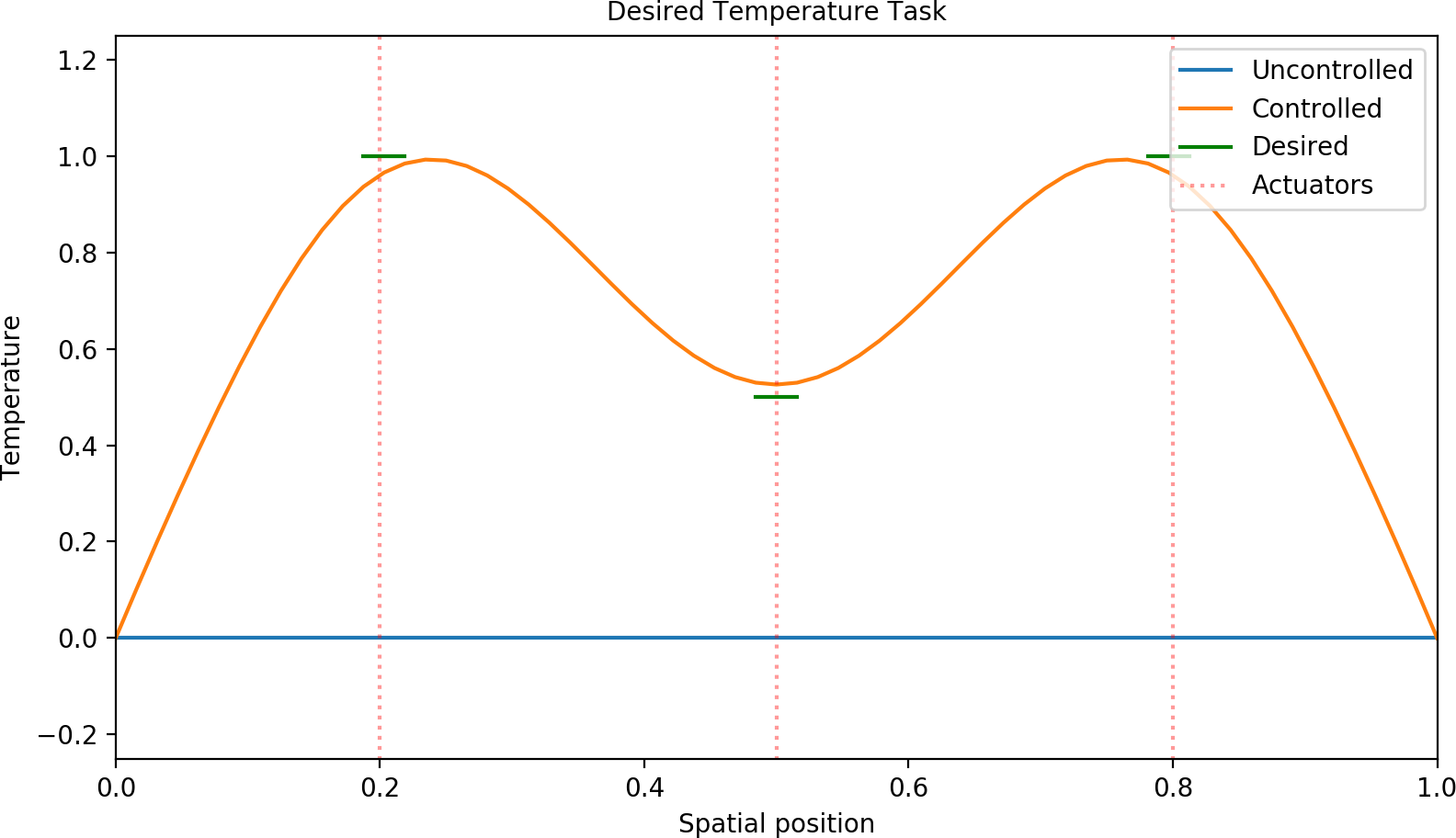}
    \end{subfigure}
    \end{multicols}
    
    \begin{subfigure}[h!]{1.0\textwidth}
    \centering
    \includegraphics[width=0.52\textwidth]{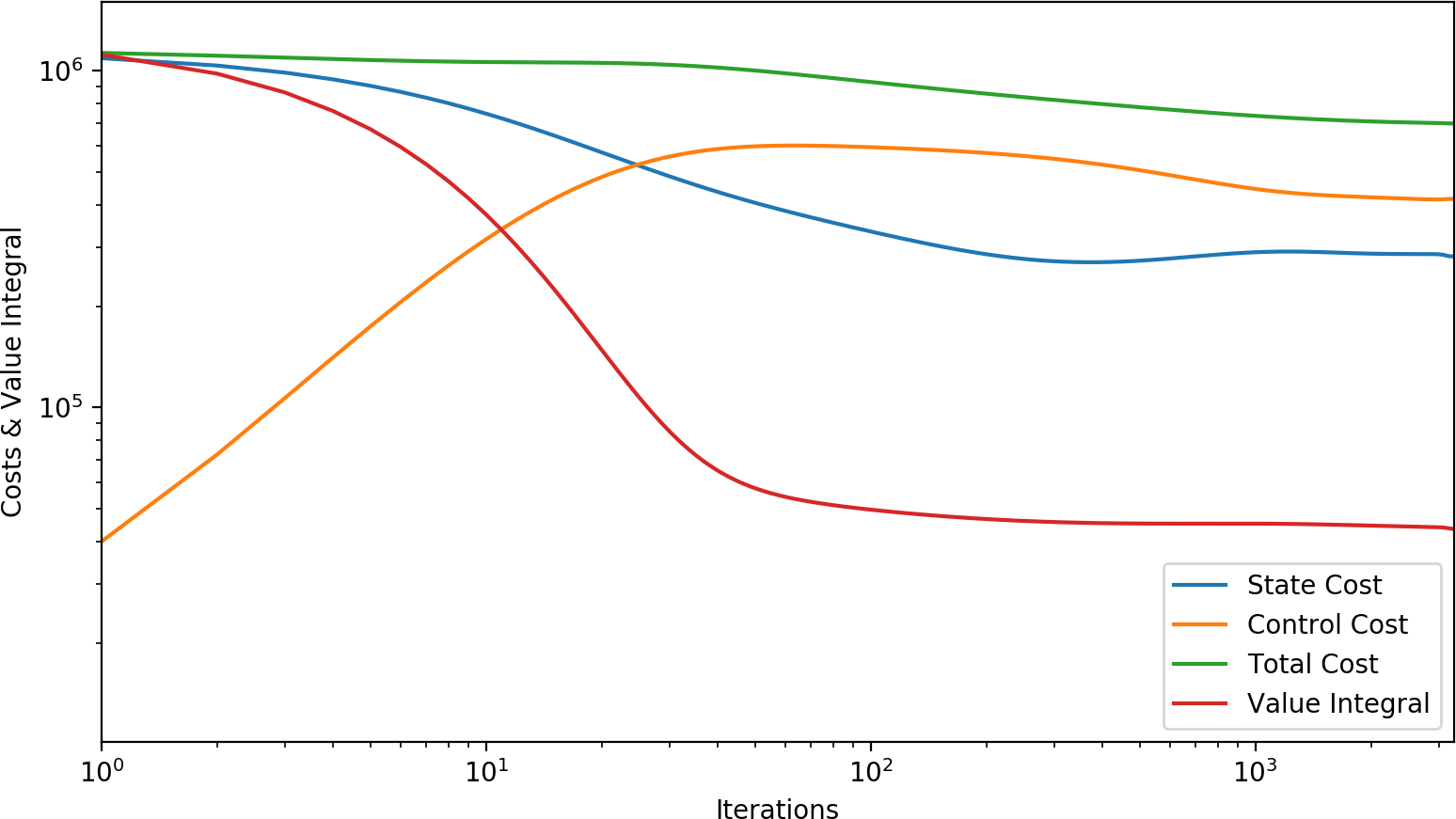}
    \end{subfigure}
    \caption{Heat Equation Temperature Reaching Task. (left) controlled contour plot where color represents temperature, (right) final time snapshot of the uncontrolled and controlled systems, (bottom) convergence plot of the heat equation temperature reaching task on a log-log scale, where the value integral depicted in red is the time integral of the value functional.}
    \label{fig:heat}
\end{figure*}}

% {\begin{figure*}[!h]
%     \centering
%     \includegraphics[width=0.75\textwidth]{Figures/1DHeat_convergence_tfpt12_scost281456_scost_w300_cropped.png}
%     \caption{Convergence plot of the heat equation temperature reaching task on a log-log scale, where the value integral depicted in red is the time integral of the value functional.}
%     \label{fig:heat_conv}
%     % \vspace{-0.7em}
% \end{figure*}}
%%%%%%%%%%%%%%HEAT FIGURES%%%%%%%%%%%%%%%%

Each experiment considered a pure quadratic cost functional of the form
\begin{align}
    J\big(t, h(t,x), U_d(t), U_b(t)\big) :=
    &\quad \Langle  h(t_f,x) - h_{\text{des}} (t_f,x), Q_f \big(h(t_f,x) - h_{\text{des}} (t_f,x)\big) \Rangle_{D_{\text{des}}} \nonumber\\
    &\quad + \int_{t_0}^{t_f} \bigg( \Langle  h(t,x) - h_{\text{des}} (t,x), Q \big(h(t,x) - h_{\text{des}} (t,x)\big) \Rangle_{D_{\text{des}}} \nonumber\\
    &\qquad\qquad\;\;  + \Langle U_d(t,x), R_d U_d(t,x)  \Rangle + \Langle U_b(t,x), R_b U_b(t,x)  \Rangle_S \bigg)\rd t,
\end{align}
where the inner product $\big\langle \cdot, \cdot \big\rangle_{D_{\text{des}}}$ is defined on the desired subregion $D_{\text{des}} \subseteq D$.

The first experiment was a temperature reaching task on the \ac{1D} Heat equation with homogeneous Dirichlet boundary conditions, given in fields representation by
\begin{equation}
    \begin{split}
    \partial_t h(t,x) &= \epsilon \partial_{xx} h(t,x) + \vm(\vx)^{\top} \, U_d(t,x),\\
    h(t,0) &= h(t,a) = 0, \\
    h(0,x) &= h_0(x),
    \end{split}
\end{equation}
where $\epsilon$ is the thermal diffusivity parameter. 
The heat equation is a pure diffusion equation, and validates the approach's ability to achieve high quality distributed control solutions in the linear \acp{PDE} regime. The \ac{STDDP} algorithm was run until convergence, and the results of which are depicted in \cref{fig:heat}. Starting from a zero initial condition, the \ac{PDE} was tasked with raising the temperature to $T=1.0$ at the outer regions, and raising the temperature to $T=0.5$ at the central region. 

The system was temporally discretized into $1200$ time steps and spatially discretized into $64$ grid points. The typical convergence behavior for the \ac{STDDP} algorithm applied to the heat equation is depicted in the bottom subfigure of \cref{fig:heat}. In this case, the weight values were $R_d = 0.4$, $Q = 300$, and $Q_f = 300$. Depicted is a log-log plot of the cost functional $J\big(t, X(t)\big)$, its state cost functional and control cost functional components, and the time integral of the value functional, which is concisely termed the value integral. The convergence behavior of the value integral demonstrates super-quadratic convergence in the first $50$ iterations.

%%%%%%%%%%%%%%%%%%%%%%%%%BURGERS%%%%%%%%%%%%%%%%%%%%%%%%%%%%%%

The second experiment was a velocity reaching task on the \ac{1D} Burgers equation with non-homogenous Dirichlet boundary conditions, given in \textit{fields representation} by
\begin{equation} \label{eq:BurgersSPDE}
\begin{split}
\partial_t h(t, x) &= - h(t,x) \partial_x h(t, x) + \epsilon \partial_{xx} h(t,x) + \vm(\vx)^{\top} \, U_d(t,x), \\
h(t,0) &= h(t,a) = 1.0,\\ 
h(0,x) &= h_0(x), 
\end{split}
\end{equation}
where the parameter $\epsilon$ is the viscosity of the medium. The Burgers equation is a nonlinear \ac{PDE}, and demonstrates the efficacy of the approach on nonlinear \acp{PDE}. Starting from a zero intiail condition, the \ac{PDE} is tasked with raising the velocity to $v=2.0$ on the outer regions, and $v=1.0$ on the central region. The Burgers equation is often used as a simplified model of fluid flow, however also has applications in describing the dynamics of swarms for robotic systems \cite{elamvazhuthi2018pde}. The \ac{STDDP} was applied to the Burgers \ac{PDE} and was run until convergence. The results are depicted in \cref{fig:burgers}.

%%%%%%%%%%%%%%BURGERS FIGURES%%%%%%%%%%%%%%%%
{\centering
\begin{figure*}[!t]
    % \subfigure[Controlled Contour]{\includegraphics[width=0.52\textwidth]{Figures/1DBurgers_contour_scost13672_cropped.png}}
    % \subfigure[Final Time Snapshot]{\hspace{.4cm}\includegraphics[width=0.453\textwidth]{Figures/1DBurgers_final_time_scost13672_cropped.png}}
    
    \begin{multicols}{2}
    \begin{subfigure}[h!]{1.0\textwidth}
    \hspace{0.0cm}\includegraphics[width=0.52\textwidth]{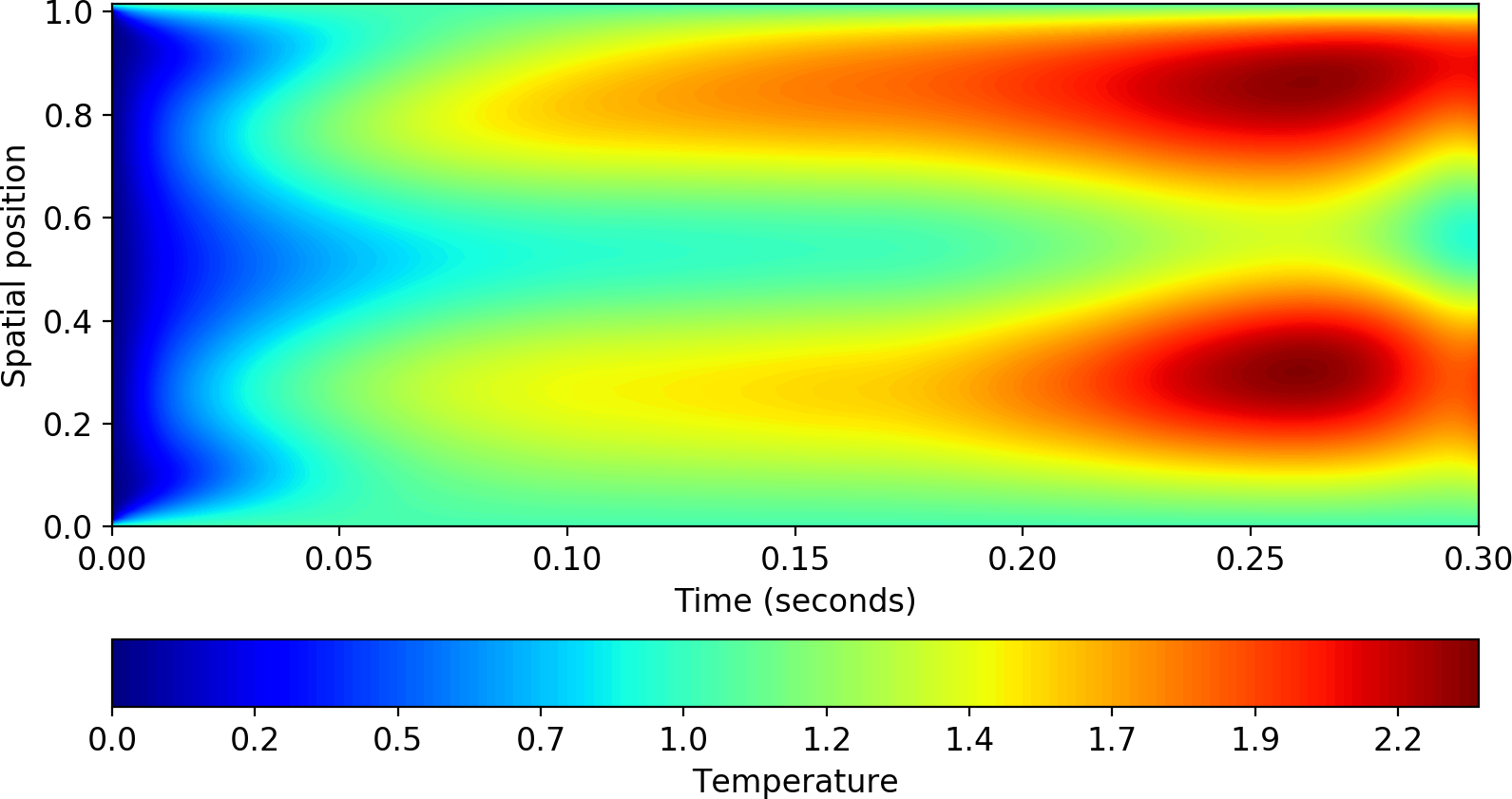}
    \end{subfigure}
    
    \begin{subfigure}[h!]{1.0\textwidth}
     \hspace{.4cm}\includegraphics[width=0.46\textwidth]{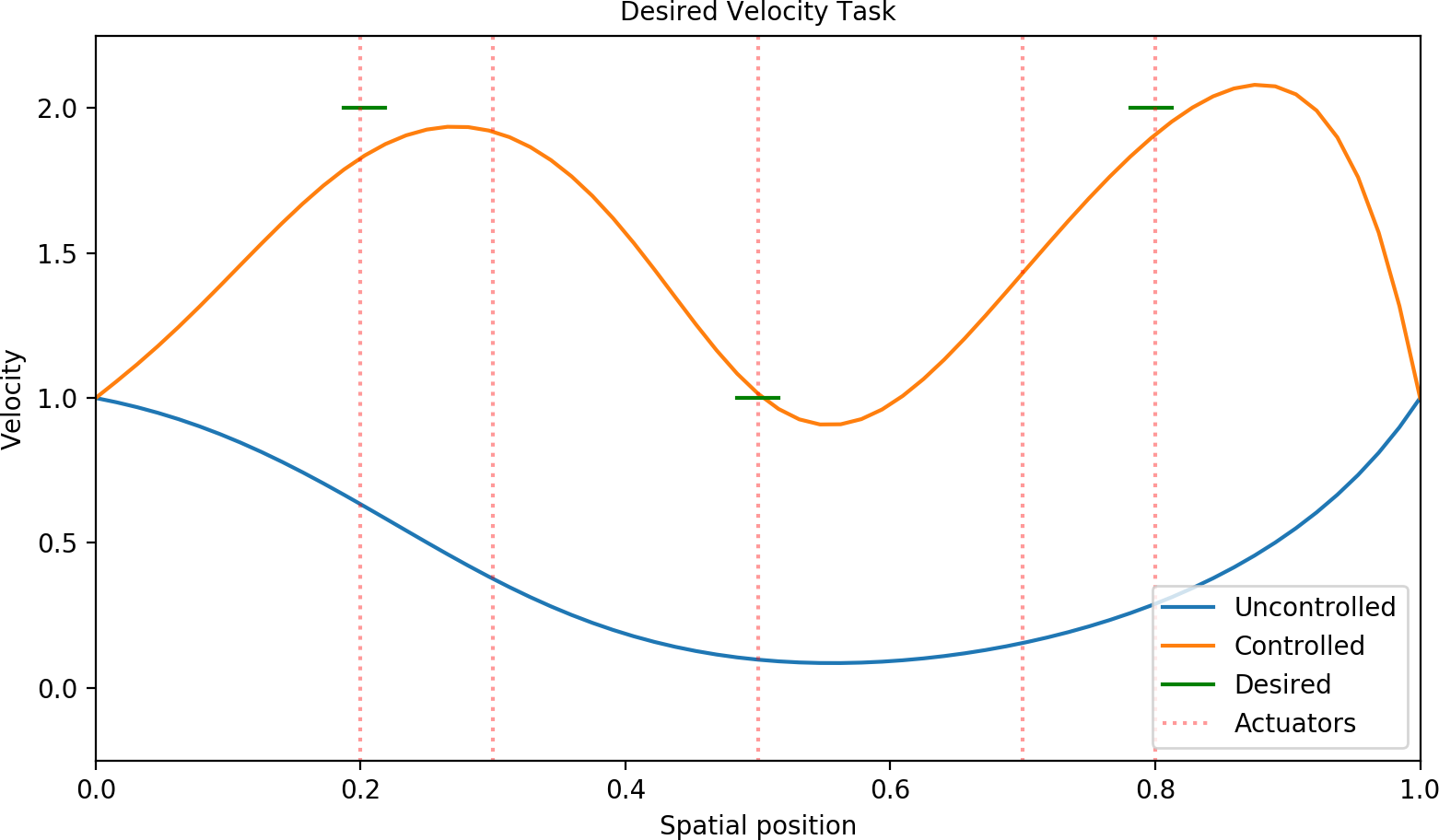}
    \end{subfigure}
    \end{multicols}
    \caption{Burgers Equation Velocity Reaching Task. (left) controlled contour plot where color represents velocity, (right) final time snapshot comparing to the uncontrolled system.}
    \label{fig:burgers}
    \vspace{-0em}
\end{figure*}}
%%%%%%%%%%%%%%BURGERS FIGURES%%%%%%%%%%%%%%%%

%%%%%%%%%%%%%%%%%%ANNEALING FIGURES%%%%%%%%%%%%%%%%%%%%%%
{\centering
\begin{figure*}[!t]
    % \subfigure[Controlled Contour]{\includegraphics[width=0.52\textwidth]{Figures/1D Burgers state cost_iter4626_costannealing1_cropped.png}}
    % \subfigure[Final Time Snapshot]{\hspace{0.4
    %  cm}\includegraphics[width=0.453\textwidth]{Figures/1D Burgers state cost_iter4626_costannealing_snapshot1_cropped.png}}
    
    \begin{multicols}{2}
    \begin{subfigure}[h!]{1.0\textwidth}
    \hspace{0.0cm}\includegraphics[width=0.52\textwidth]{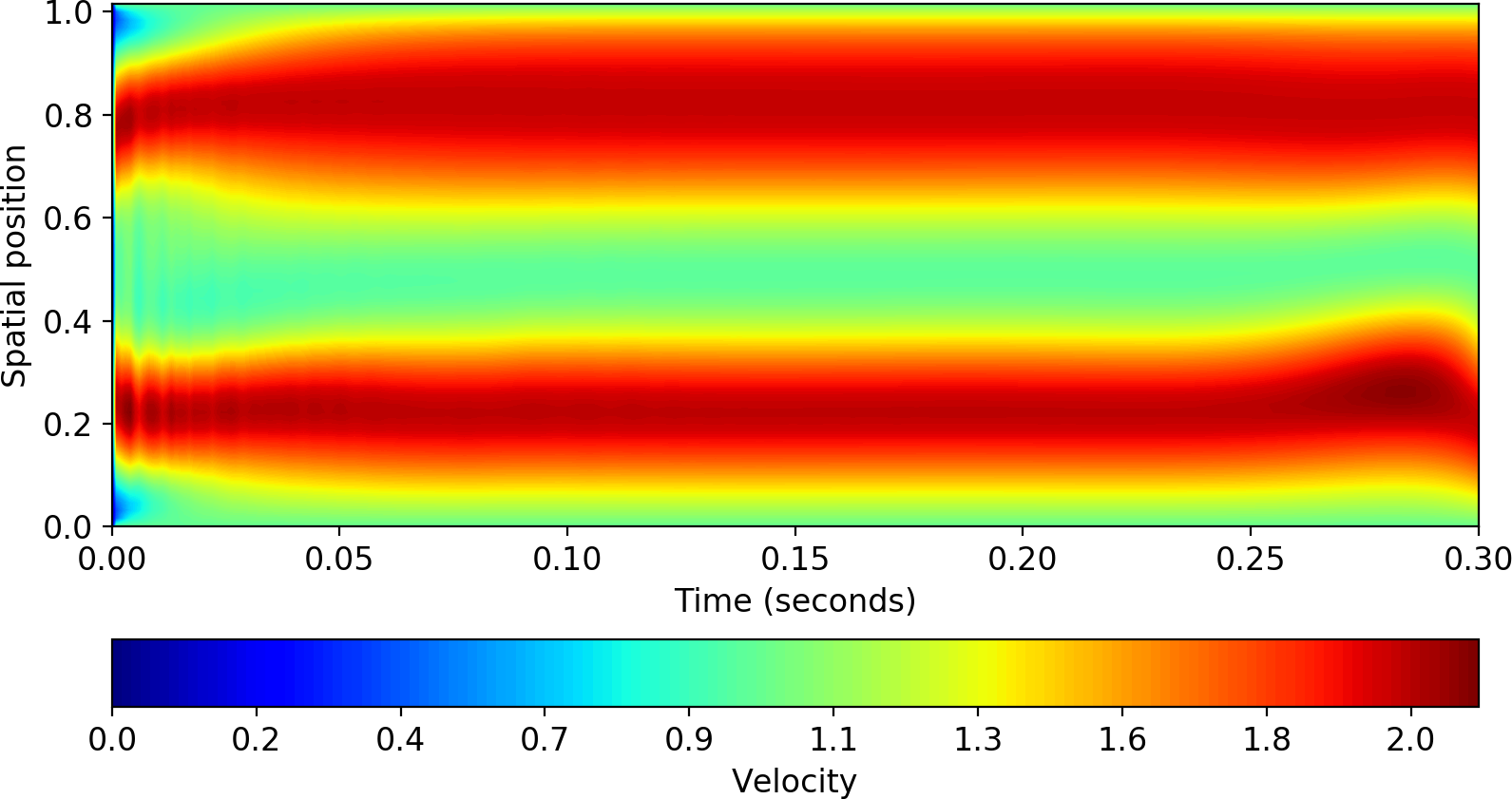}
    \end{subfigure}
    
    \begin{subfigure}[h!]{1.0\textwidth}
     \hspace{0.4cm}\includegraphics[width=0.46\textwidth]{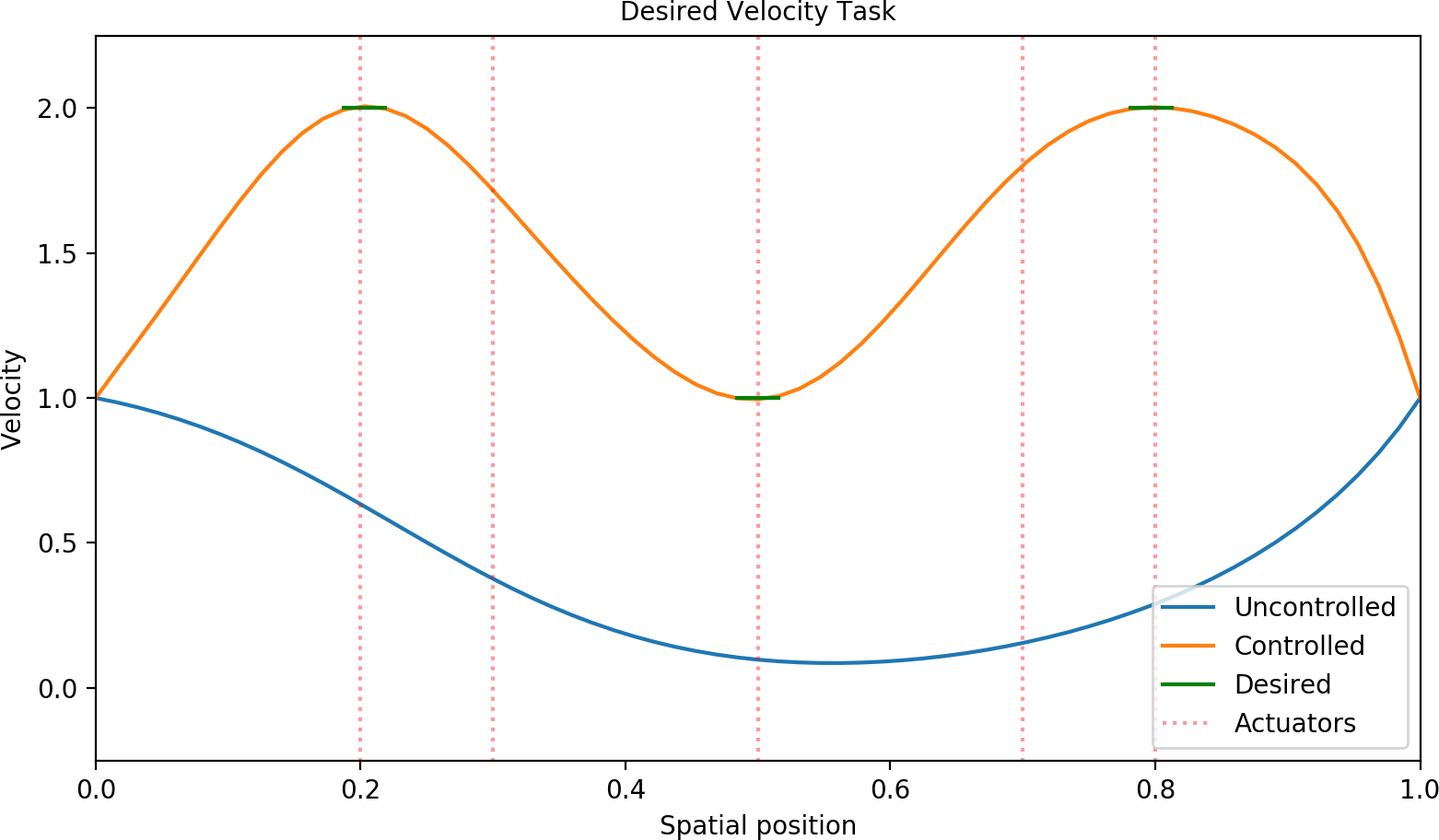}
    \end{subfigure}
    \end{multicols}
    \caption{Burgers Equation Velocity Reaching Task with Simulated Annealing. (left) controlled contour plot where color represents velocity, (right) final time snapshot comparing the optimized solution to the uncontrolled system.}
    \label{fig:burgers_annealing}
\end{figure*}}
%%%%%%%%%%%%%%%%%%ANNEALING FIGURES%%%%%%%%%%%%%%%%%%%%%%

The nonlinear advection present in the Burgers equation produces an apparent rightward motion that builds over the spatial domain to create an apparent wavefront towards the right endpoint. The system is provided with $5$ actuators, and must overcome this nonlinearity in order to minimize the state cost. Despite the added actuators, the task remains severely under-actuated. In this case, the weight values were $R_d = 0.4$, $Q=30$, and $Q_f = 30$. As depicted, the provided values of state and control cost weighting provide a balancing between the state and control performance metrics.

In both of the experiments, the various discretization schemes describbed in \cref{sec:discretization} were tested, namely the explicit Euler discretization, a Runge-Kutta 2-point discretization, and the semi-implicit. The authors report that while the semi-implicit method had slightly lower sensitivity to the time-step increment $\Delta t$ compared to the explicit Euler and Runge-Kutta methods, the large matrix inversion caused dramatically slower per-iteration run-time. The Runge-Kutta method had higher accuracy than the Euler method, but required super-sampling (i.e. sampling the midpoint of a time-increment) thus doubling the total time steps on forward and backward passes. The explicit Euler discretization had the fastest per-iteration run time at about $0.4$ seconds per iteration, and was stabilized using regularization methods, akin to \cite{tassa2012synthesis}.
% \textbf{maybe make a table of per-iteration run-time, dt sensitivity, and accuracy}

Common to finite and infinite dimensional \ac{DDP} methods are parameter sensitivities that may limit choice of the control cost weighting and the state cost weightings. When these limits arise, they are typically due to the numerically stiff and sensitive dynamics found in the backwards Ricatti equation \cref{eq:second_riccati}, and present a limitation in the ability of \ac{DDP} approaches to use arbitrary ratios of state performance and control effort. This can be especially limiting in systems with under-actuation as control signals can often be much larger for task completion, thus requiring a larger ratio between state cost weight and control cost weight. Without a "warm start", the operational initialization window for control weights may limit the use of an arbitrary desired set of parameters, thus changing the task specifications to meet numerical requirements.

In \cref{fig:burgers_annealing}, we demonstrate that this can be overcome with a simple simulated annealing scheme. In this simulated experiment, the simulated annealing scheme was adopted in order to reach an arbitrarily large weight ratio $W_d := Q/R_d = 4.8\times10^6$ starting from a nominal weight ratio of $W_d = 25$. This approach allows one to arbitrarily choose the relative importance of state performance and control effort. Depicted is a contour plot that demonstrates that the desired regions are quickly reached, and the system remains at the desired region for the duration of the simulation. Also depicted is a final time snapshot with dramatically smaller deviation from the desired region as compared to the solution in \cref{fig:burgers}, albeit at the expense of larger control effort.

%===============================================================================

\section{Discussion \& Conclusion}

We address the optimal control on nonlinear spatio-temporal systems through the lens of the Bellman principle of optimality, and develop the \ac{STDDP} framework. We demonstrate that the resulting forward-backward system of equations can recover standard results, including the \ac{LQR} solution for linear \acp{PDE} and the \ac{DDP} solution for finite nonlinear \acp{ODE}. We analyze the convergence behavior and emerge with provable global convergence of the resulting forward-backward system. We discuss and develop discretization schemes for the backward second derivative of the value functional, and implement the resulting algorithm on a linear \ac{PDE} system and a nonlinear \ac{PDE} system.

The numerical results demonstrate the utility of the \ac{STDDP} framework. It has the capability of obtaining high quality control solutions in the linear and nonlinear regime for spatio-temporal \ac{PDE} systems. It has flexibility with respect to discretization schemes due to the optimize-then-discretize approach. It exhibits computationally efficiency for \ac{1D} \acp{PDE} with a typical $~0.5$ second time-per-iteration without any parallelization. 

Overall, the results presented in this manuscript are encouraging to the authors for future work on extending the approach to 2D and 3D problem spaces. Such scaling will result in large tensors, however one can leverage the sparsity inherent in \ac{PDE} discretizations and utilize common tensor decompositions such as the tensor train decomposition \cite{oseledets2011tensor} for a dramatic computational speed-up. Other future directions include extensions to the case of a system with additive Gaussian noise, second order expansions of the dynamics, and novel methods to handle the sensitivities that arise in the discretization of the backward process.

%===============================================================================

\section*{Acknowledgements}
This work was supported by the Army Research Office contract W911NF2010151. Ethan N. Evans was supported by the SMART scholarship.

\bibliographystyle{ieeetran}
\balance
\bibliography{References}

\begin{thebibliography}{10}
\providecommand{\url}[1]{#1}
\csname url@samestyle\endcsname
\providecommand{\newblock}{\relax}
\providecommand{\bibinfo}[2]{#2}
\providecommand{\BIBentrySTDinterwordspacing}{\spaceskip=0pt\relax}
\providecommand{\BIBentryALTinterwordstretchfactor}{4}
\providecommand{\BIBentryALTinterwordspacing}{\spaceskip=\fontdimen2\font plus
\BIBentryALTinterwordstretchfactor\fontdimen3\font minus
  \fontdimen4\font\relax}
\providecommand{\BIBforeignlanguage}[2]{{%
\expandafter\ifx\csname l@#1\endcsname\relax
\typeout{** WARNING: IEEEtran.bst: No hyphenation pattern has been}%
\typeout{** loaded for the language `#1'. Using the pattern for}%
\typeout{** the default language instead.}%
\else
\language=\csname l@#1\endcsname
\fi
#2}}
\providecommand{\BIBdecl}{\relax}
\BIBdecl

\bibitem{lord_powell_shardlow_2014}
G.~J. Lord, C.~E. Powell, and T.~Shardlow, \emph{An Introduction to
  Computational Stochastic PDEs}, ser. Cambridge Texts in Applied
  Mathematics.\hskip 1em plus 0.5em minus 0.4em\relax Cambridge University
  Press, 2014.

\bibitem{gomes2017controlling}
S.~N. Gomes, S.~Kalliadasis, D.~T. Papageorgiou, G.~A. Pavliotis, and
  M.~Pradas, ``Controlling roughening processes in the stochastic
  kuramoto--sivashinsky equation,'' \emph{Physica D: Nonlinear Phenomena}, vol.
  348, pp. 33--43, 2017.

\bibitem{rabault2019artificial}
J.~Rabault, M.~Kuchta, A.~Jensen, U.~R{\'e}glade, and N.~Cerardi, ``Artificial
  neural networks trained through deep reinforcement learning discover control
  strategies for active flow control,'' \emph{Journal of Fluid Mechanics}, vol.
  865, pp. 281--302, 2019.

\bibitem{bieker2019deep}
K.~Bieker, S.~Peitz, S.~L. Brunton, J.~N. Kutz, and M.~Dellnitz, ``Deep model
  predictive control with online learning for complex physical systems,''
  \emph{arXiv preprint arXiv:1905.10094}, 2019.

\bibitem{mohan2018deep}
A.~T. Mohan and D.~V. Gaitonde, ``A deep learning based approach to reduced
  order modeling for turbulent flow control using lstm neural networks,''
  \emph{arXiv preprint arXiv:1804.09269}, 2018.

\bibitem{nair2019cluster}
A.~G. Nair, C.-A. Yeh, E.~Kaiser, B.~R. Noack, S.~L. Brunton, and K.~Taira,
  ``Cluster-based feedback control of turbulent post-stall separated flows,''
  \emph{Journal of Fluid Mechanics}, vol. 875, pp. 345--375, 2019.

\bibitem{satheeshbabu2019open}
S.~Satheeshbabu, N.~K. Uppalapati, G.~Chowdhary, and G.~Krishnan, ``Open loop
  position control of soft continuum arm using deep reinforcement learning,''
  in \emph{2019 International Conference on Robotics and Automation
  (ICRA)}.\hskip 1em plus 0.5em minus 0.4em\relax IEEE, 2019, pp. 5133--5139.

\bibitem{spielberg2019learning}
A.~Spielberg, A.~Zhao, Y.~Hu, T.~Du, W.~Matusik, and D.~Rus,
  ``Learning-in-the-loop optimization: End-to-end control and co-design of soft
  robots through learned deep latent representations,'' \emph{Advances in
  Neural Information Processing Systems}, vol.~32, pp. 8284--8294, 2019.

\bibitem{farahmand2017deep}
A.-m. Farahmand, S.~Nabi, and D.~N. Nikovski, ``Deep reinforcement learning for
  partial differential equation control,'' in \emph{2017 American Control
  Conference (ACC)}.\hskip 1em plus 0.5em minus 0.4em\relax IEEE, 2017, pp.
  3120--3127.

\bibitem{morton2018deep}
J.~Morton, A.~Jameson, M.~J. Kochenderfer, and F.~Witherden, ``Deep dynamical
  modeling and control of unsteady fluid flows,'' in \emph{Advances in Neural
  Information Processing Systems}, 2018, pp. 9258--9268.

\bibitem{lasiecka2000control}
I.~Lasiecka and R.~Triggiani, \emph{Control theory for partial differential
  equations: continuous and approximation theories}.\hskip 1em plus 0.5em minus
  0.4em\relax Cambridge University Press Cambridge, 2000, vol.~1.

\bibitem{troltzsch2010optimal}
F.~Tr{\"o}ltzsch, \emph{Optimal control of partial differential equations:
  theory, methods, and applications}.\hskip 1em plus 0.5em minus 0.4em\relax
  American Mathematical Soc., 2010, vol. 112.

\bibitem{sumin2009first}
M.~I. Sumin, ``The first variation and pontryagin’s maximum principle in
  optimal control for partial differential equations,'' \emph{Computational
  Mathematics and Mathematical Physics}, vol.~49, no.~6, pp. 958--978, 2009.

\bibitem{yong1992pontryagin}
J.~M. Yong, ``Pontryagin maximum principle for semilinear second order elliptic
  partial differential equations and variational inequalities with state
  constraints,'' \emph{Differential and Integral Equations}, vol.~5, no.~6, pp.
  1307--1334, 1992.

\bibitem{tassa2014control}
Y.~Tassa, N.~Mansard, and E.~Todorov, ``Control-limited differential dynamic
  programming,'' in \emph{IEEE International Conference on Robotics and
  Automation (ICRA)}, May 2014, pp. 1168--1175.

\bibitem{aoyama2020constrained}
Y.~Aoyama, G.~Boutselis, A.~Patel, and E.~A. Theodorou, ``Constrained
  differential dynamic programming revisited,'' \emph{arXiv preprint
  arXiv:2005.00985}, 2020.

\bibitem{tassa2007receding}
Y.~Tassa, T.~Erez, and W.~D. Smart, ``Receding horizon differential dynamic
  programming.'' in \emph{NIPS}.\hskip 1em plus 0.5em minus 0.4em\relax
  Citeseer, 2007, pp. 1465--1472.

\bibitem{pan2014probabilistic}
Y.~Pan and E.~Theodorou, ``Probabilistic differential dynamic programming,''
  \emph{Advances in Neural Information Processing Systems}, vol.~27, pp.
  1907--1915, 2014.

\bibitem{pan2018efficient}
Y.~Pan, G.~I. Boutselis, and E.~A. Theodorou, ``Efficient reinforcement
  learning via probabilistic trajectory optimization,'' \emph{IEEE transactions
  on neural networks and learning systems}, vol.~29, no.~11, pp. 5459--5474,
  2018.

\bibitem{sun2015game}
W.~Sun, E.~A. Theodorou, and P.~Tsiotras, ``Game theoretic continuous time
  differential dynamic programming,'' in \emph{2015 American Control Conference
  (ACC)}.\hskip 1em plus 0.5em minus 0.4em\relax IEEE, 2015, pp. 5593--5598.

\bibitem{boutselis2019numerical}
G.~I. Boutselis, Y.~Pan, and E.~A. Theodorou, ``Numerical trajectory
  optimization for stochastic mechanical systems,'' \emph{SIAM Journal on
  scientific computing}, vol.~41, no.~4, pp. A2065--A2087, 2019.

\bibitem{boutselis2016stochastic}
G.~I. Boutselis, G.~De~La~Torre, and E.~A. Theodorou, ``Stochastic optimal
  control using polynomial chaos variational integrators,'' in \emph{2016
  American Control Conference (ACC)}.\hskip 1em plus 0.5em minus 0.4em\relax
  IEEE, 2016, pp. 6586--6591.

\bibitem{tzafestas1969differential}
S.~Tzafestas and J.~Nightingale, ``Differential dynamic-programming approach to
  optimal nonlinear distributed-parameter control systems,'' in
  \emph{Proceedings of the Institution of Electrical Engineers}, vol. 116,
  no.~6.\hskip 1em plus 0.5em minus 0.4em\relax IET, 1969, pp. 1079--1084.

\bibitem{sakawa1972matrix}
Y.~Sakawa, ``A matrix green's formula and optimal control of linear
  distributed-parameter systems,'' \emph{Journal of Optimization Theory and
  Applications}, vol.~10, no.~5, pp. 290--299, 1972.

\bibitem{jacobson1970}
D.~H. Jacobson and D.~Q. Mayne, \emph{Differential dynamic programming}.\hskip
  1em plus 0.5em minus 0.4em\relax New York: American Elsevier Pub. Co., 1970.

\bibitem{evans1997partial}
L.~C. Evans, ``Partial differential equations and monge-kantorovich mass
  transfer,'' \emph{Current developments in mathematics}, vol. 1997, no.~1, pp.
  65--126, 1997.

\bibitem{volterra1959theory}
V.~Volterra, \emph{Theory of functionals and of integral and
  integro-differential equations}.\hskip 1em plus 0.5em minus 0.4em\relax
  Dover, 1959.

\bibitem{shoemaker1990proof}
C.~Shoemaker and L.~Liao, ``Proof of the quadratic convergence of differential
  dynamic programming,'' Cornell University Operations Research and Industrial
  Engineering, Tech. Rep., 1990.

\bibitem{yakowitz1984computational}
S.~Yakowitz and B.~Rutherford, ``Computational aspects of discrete-time optimal
  control,'' \emph{Applied Mathematics and Computation}, vol.~15, no.~1, pp.
  29--45, 1984.

\bibitem{sun2014continuous}
W.~Sun, E.~A. Theodorou, and P.~Tsiotras, ``Continuous-time differential
  dynamic programming with terminal constraints,'' in \emph{2014 IEEE Symposium
  on Adaptive Dynamic Programming and Reinforcement Learning (ADPRL)}.\hskip
  1em plus 0.5em minus 0.4em\relax IEEE, 2014, pp. 1--6.

\bibitem{tzafestas1968optimal}
S.~Tzafestas and J.~Nightingale, ``Optimal control of a class of linear
  stochastic distributed-parameter systems,'' in \emph{Proceedings of the
  Institution of Electrical Engineers}, vol. 115, no.~8.\hskip 1em plus 0.5em
  minus 0.4em\relax IET, 1968, pp. 1213--1220.

\bibitem{friedman2008partial}
A.~Friedman, \emph{Partial differential equations of parabolic type}.\hskip 1em
  plus 0.5em minus 0.4em\relax Courier Dover Publications, 2008.

\bibitem{pantoja1983algorithms}
J.~D.~O. Pantoja, ``Algorithms for constrained optimization problems,''
  \emph{Differential Dynamic Programming and Newton's Method. International
  Journal of Control}, vol.~47, pp. 1539--1553, 1983.

\bibitem{murray1984differential}
D.~Murray and S.~Yakowitz, ``Differential dynamic programming and newton's
  method for discrete optimal control problems,'' \emph{Journal of Optimization
  Theory and Applications}, vol.~43, no.~3, pp. 395--414, 1984.

\bibitem{boutselis2018differential}
G.~I. Boutselis and E.~Theodorou, ``Differential dynamic programming on lie
  groups: Derivation, convergence analysis and numerical results,'' \emph{arXiv
  preprint arXiv:1809.07883}, 2018.

\bibitem{tassa2012synthesis}
Y.~Tassa, T.~Erez, and E.~Todorov, ``Synthesis and stabilization of complex
  behaviors through online trajectory optimization,'' in \emph{IEEE/RSJ
  International Conference on Intelligent Robots and Systems (IROS)}, Oct 2012,
  pp. 4906--4913.

\bibitem{butusov2016semi}
D.~Butusov, T.~Karimov, and V.~Ostrovskii, ``Semi-implicit ode solver for
  matrix riccati equation,'' in \emph{2016 IEEE NW Russia Young Researchers in
  Electrical and Electronic Engineering Conference (EIConRusNW)}.\hskip 1em
  plus 0.5em minus 0.4em\relax IEEE, 2016, pp. 168--172.

\bibitem{demmel1999supernodal}
J.~W. Demmel, S.~C. Eisenstat, J.~R. Gilbert, X.~S. Li, and J.~W. Liu, ``A
  supernodal approach to sparse partial pivoting,'' \emph{SIAM Journal on
  Matrix Analysis and Applications}, vol.~20, no.~3, pp. 720--755, 1999.

\bibitem{elamvazhuthi2018pde}
K.~Elamvazhuthi, H.~Kuiper, and S.~Berman, ``Pde-based optimization for
  stochastic mapping and coverage strategies using robotic ensembles,''
  \emph{Automatica}, vol.~95, pp. 356--367, 2018.

\bibitem{oseledets2011tensor}
I.~V. Oseledets, ``Tensor-train decomposition,'' \emph{SIAM Journal on
  Scientific Computing}, vol.~33, no.~5, pp. 2295--2317, 2011.

\end{thebibliography}

%===============================================================================

\newpage
\section*{Supplementary Information}
\beginsupplement

\section{Proof of Lemma 9.1}\label{sec:proof_convergence_lemma}

\begin{proof}
We start with the total derivative for the cost functional
\begin{align}
    \frac{\rd J^i\big(t, X(t), U(t)\big) }{ \rd U^i_{t_0:t_f}} &=  \frac{\partial J^i\big(t, X(t),  U(t)\big) }{ \partial U^i_{t_0:t_f}} + \frac{\partial J^i\big(t, X(t),  U(t)\big) }{ \partial X^i_{t_0:t_f}}  \frac{\partial X^i_{t_0:t_f} }{ \partial U^i_{t_0:t_f}} \\
    &= \Langle L_{U,t_0:t_f}, \mathds{1}_{t_0:t_f} \Rangle_T + \Langle L_{X,t_0:t_f},  \frac{\partial X^i_{t_0:t_f} }{ \partial U^i_{t_0:t_f}} \Rangle_T + \Langle \phi_X^\top\big(t_f, X(t_f)\big), \frac{\partial X^i_{t_f} }{ \partial U^i_{t_0:t_f}} \Rangle
\end{align}
The state trajectory $X_{t_0:t_f}$ is due to the approximate state evolution given by \cref{eq:F_expand,eq:N_expand}, which has linear affine form with solution
\begin{equation}
    X(t) = \Phi(t,t_0)X_0 + \int_{t_0}^{t} \Phi(t,s) F_U(s) U(s) \rd s
\end{equation}
Thus, one has
\begin{align}
    \frac{\rd J^i\big(t, X(t), U(t)\big) }{ \rd U^i_{t_0:t_f}} &= \Langle L_{U,t_0:t_f}, \mathds{1}_{t_0:t_f} \Rangle_T + \Langle \Phi^\top(t,t_0:t_f) L_{X,t_0:t_f},  F_{U,t_0:t_f} \Rangle_T \nonumber \\
    &\quad + \Langle \Phi^\top(T,t_0:t_f) \phi_X^\top\big(t_f, X(t_f)\big),  F_{U,t_0:t_f} \Rangle_T
\end{align}
Now, due to the terminal condition on $\psi$ given in \cref{eq:psi_BC}, one has
\begin{align}
    \frac{\rd J^i\big(t, X(t), U(t)\big) }{ \rd U^i_{t_0:t_f}} &= \Langle L_{U,t_0:t_f}, \mathds{1}_{t_0:t_f} \Rangle_T + \Langle \Phi^\top(t,t_0:t_f) L_{X,t_0:t_f},  F_{U,t_0:t_f} \Rangle_T \nonumber \\
    &\quad + \Langle \Phi^\top(T,t_0:t_f) \psi\big(t_f, X(t_f), U(t_f)\big),  F_{U,t_0:t_f} \Rangle_T.
\end{align}
Finally, due to the backward recursion over the trajectory given by \cref{eq:psi}, one has
\begin{align}
    \frac{\rd J^i\big(t, X(t), U(t)\big) }{ \rd U^i_{t_0:t_f}} = \Langle L_{U,t_0:t_f}, \mathds{1}_{t_0:t_f} \Rangle_T + \Langle \psi_{t_0:t_f},  F_{U,t_0:t_f} \Rangle_T,
\end{align}
which concludes the proof.
\end{proof}

\section{Proof of Corollary 9.1}\label{sec:proof_convergence_corollary}
\begin{proof}
For simplicity of notation, we will use the shorthand $J^i := J^i\big(t,X^i(t),U^i(t)\big)$. Theorem \ref{thm:quadratic} and \cref{assumption_compactness} together imply that $\exists \gamma \in [0,1)$ such that the change in the cost over iterations $\Delta J^i := J^i - J^{i-1} < 0$, $\forall i \in \Nb_+$. The cost functional $J(\cdot,\cdot,\cdot)$ is continuous in its arguments since it is differentiable by assumption, therefore it is also continuous with respect to iterations $i$. Thus $\lim_{i\rightarrow \infty} \Delta J^i = 0$ and $\exists$ a pair $(X^*_{t_0:t_f}, U^*_{t_0:t_f})$ such that $\lim_{i \rightarrow \infty} J\big(t,X^i(t), U^i(t) \big) = J\big(t, X^*(t), U^*(t) \big) =: J^*$.

Next, $\lim_{i\rightarrow \infty} \Delta J^i = 0$ together with \cref{eq:conv_thm} and the positive definiteness of $M_{t_0:t_f}$ imply that \\
$\lim_{i\rightarrow \infty} Q^i_{U,t_0:t_f} = 0_{t_0:t_f}$. By this notation we mean that $\lim_{i\rightarrow \infty} Q^i_U\big(t, X(t), U(t)\big) =  0$ $\forall t \in T$. Recall also that $\delX_0^i$. We seek the intermediate result that $\delX^*_{t_0:t_f} := \lim_{i\rightarrow \infty} X^i_{t_0:t_f} = 0_{t_0:t_f}$. This can be observed in the coupled solutions of $\delX(t)$ and $\delU(t)$, but is more clear by inspection of the closed-loop variation dynamics, which have generalized form
\begin{align*}
    \frac{\rd \delX(t)}{\rd t} &= F_X^\top\!\!(t) \delX(t) - F_U^\top\!\!(t)Q_{UU}^{-1}(t)\Big(Q_U\!(t) + Q_{UX}\!(t)\delX(t) \Big) \\
    &= \Big(F_X^\top(t) - F_U^\top(t) Q_{UU}^{-1}(t) Q_{UX}(t)\Big)\delX(t) - F_U^\top(t)Q_{UU}^{-1}(t)Q_U(t)
\end{align*}
with solution of the form
\begin{align*}
    \delX(t) = \tilde{\Phi}(t,t_0)\delX_0 - \int_{t_0}^t \tilde{\Phi}(t,s)F_U^\top(s)Q_{UU}^{-1}(s)Q_U(s) \rd s,
\end{align*}
where $\tilde{\Phi}(\cdot,\cdot)$ is a contractive linear semigroup. Thus, since $\delX_0 = 0$, and $\lim_{i \rightarrow \infty}Q^i_{U,t_0:t_f} = 0$, we have that $\delX^*_{t_0:t_f} = 0_{t_0:t_f}$ and thus $\delU^*_{t_0:t_f} := \lim_{i\rightarrow \infty} \delU^i_{t_0:t_f} = 0_{t_0:t_f}$, which implies that $\lim_{i\rightarrow \infty}U^i_{t_0:t_f} = U^*_{t_0:t_f}$.

Finally, we must show that the converged control trajectory $U^*_{t_0:t_f}$ is stationary. To show this, consider the rate of change of the cost functional with respect to control over iterations in the limit, namely
\begin{align*}
\limiinfty \frac{\rd J^i}{\rd U^i_{t_0:t_f}} &= \limiinfty\Langle Q_{U,t_0:t_f}^i, \mathds{1}_{t_0:t_f} \Rangle_T - \limiinfty \Langle F_{U,t_0:t_f}^i, V_{X,t_0:t_f}^i - \psi_{t_0:t_f} \Rangle_T
\end{align*}
Since we already showed that $\limiinfty Q_{U,t_0:t_f}^i = 0$, one can easily apply the dominated convergence theorem to show that the first term is a zero trajectory. We must only prove that the second term is also a zero trajectory. By \cref{D_prop}, we have
\begin{align*}
    \limiinfty \Big(V_{X}^i(t) - \psi(t)\Big) &= \limiinfty \int_T^t \Phi^\top(t,\tau) Q_{UX}^{i^\top}(\tau) Q_{UU}^{i^{-1}}Q_U^i(\tau) \rd \tau \\
    &=  \int_T^t \limiinfty \Phi^\top(t,\tau) Q_{UX}^{i^\top}(\tau) Q_{UU}^{i^{-1}}Q_U^i(\tau) \rd \tau \\
    &=0
\end{align*}
where we have again applied a properly formulated dominated convergence argument due to boundedness of  $V_{XX,t_0:t_f}$ and $X_{t_0:t_f}$ $\forall i \in \Nb_+$ by assumption, and due to $Q_{U,t_0:t_f}$ being a decreasing function over iterations. The limit is again zero since $\limiinfty Q_{U,t_0:t_f}^i = 0_{t_0:t_f}$. Thus $\limiinfty \frac{\rd J^i}{\rd U^i_{t_0:t_f}}=0_{t_0:t_f}$, and the converged trajectory $U_{t_0:t_f}^*$ is indeed stationary, which concludes the proof.
\end{proof}

\end{document}